\documentclass[journal]{IEEEtran}
\usepackage{amsmath,amssymb,amsthm}
\usepackage{graphicx}
\usepackage{cite}
\usepackage{cancel}
\usepackage[normalem]{ulem}
\usepackage{subfig}
\usepackage{float}
\usepackage{mathrsfs}
\usepackage{euscript}
\usepackage{color}
\usepackage[dvipsnames]{xcolor}

\newcommand{\Cset}{\mathcal{C}}
\newcommand{\Lone}{\mathcal{L}_1}
\newcommand{\Linf}{\mathcal{L}_{\infty}}

\newcommand{\Lonenorm}[1]{\ensuremath{\left\Vert{#1}\right\Vert_{\mathcal{L}_1}}}
\newcommand{\Linfnorm}[1]{\ensuremath{\left\Vert{#1}\right\Vert_{\mathcal{L}_{\infty}}}}
\newcommand{\Linfnormt}[2]{\ensuremath{     \left\Vert{#1}_{#2}\right\Vert_{\mathcal{L}_{\infty}}}}

\newcommand{\Linfnormtb}[3]{\ensuremath{{\left\Vert{#1}\right\Vert}_{\mathcal{L}_{\infty}{[{#2},{#3}]}}}}

\newcommand{\twonorm}[1]{\ensuremath{\Vert{#1}\Vert_{2}}}
\newcommand{\norm}[1]{\ensuremath{\Vert{#1}\Vert}}
\newcommand{\set}[1]{\{{#1}\}}

\newcommand{\Proj}[2]{\ensuremath{{\rm Proj}({#1},{#2})}}

\newcommand{\ytilde}{\tilde{y}}

\newcommand{\vtilde}{\tilde{v}}
\newcommand{\xref}{{x_{ref}}}
\newcommand{\uref}{{u_{ref}}}
\newcommand{\yref}{{y_{ref}}}

\newcommand{\txref}{{\tilde{x}_{ref}}}
\newcommand{\tyref}{{\tilde{y}_{ref}}}
\newcommand{\turef}{{\tilde{u}_{ref}}}

\newcommand{\tnref}{{\tilde{\eta}_{ref}}}
\newcommand{\txrefsub}[1]{{\tilde{x}_{ {ref}_{#1} } }}

\newcommand{\turefsub}[1]{{\tilde{u}_{ {ref}_{#1} } }}

\newcommand{\IC}{{\mathbb{C}}}

\newcommand{\II}{{\mathbb{I}}}

\newcommand{\IR}{{\mathbb{R}}}

\newcommand{\Ccal}{{\mathcal{C}}}

\newtheorem{thm}{Theorem}

\newtheorem{lem}{Lemma}
\newtheorem{defn}{Definition}
\newtheorem{cor}{Corollary}
\newtheorem{asm}{Assumption}
\newtheorem{rem}{Remark}
\newtheorem{prob}{Problem}

\graphicspath{./figures/}
\begin{document}

\title{ $\Lone$ Adaptive Output Feedback for Non-square\\  Systems  with Arbitrary Relative Degree }

 \author{
 	\thanks{This work was supported by AFOSR and NASA.}
 	Hanmin Lee\thanks{Hanmin Lee is a senior researcher at the Agency of Defense Development (ADD), Daejeon, Republic of Korea, (e-mail: orangeus170@gmail.com).},
 	Venanzio Cichella\thanks{Venanzio Cichella is a professor at the Department of Mechanical Engineering, University of Iowa, Iowa City, IA 52242
 		(e-mail: venanzio-cichella@uiowa.edu).},
 	Naira Hovakimyan\thanks{Naira Hovakimyan is a professor at the Department of Mechanical Science and Engineering, University of Illinois at Urbana-Champaign, Urbana, IL 61801 (e-mail: nhovakim@illinois.edu).}
}

\maketitle

\begin{abstract}
This paper considers the problem of output-feedback control for non-square multi-input multi-output systems with arbitrary relative degree. The proposed controller, based on the $\Lone$ adaptive control architecture, is designed using the right interactor matrix and a suitably defined projection matrix. A state-output predictor, a low-pass filter, and adaptive laws are introduced that achieve output tracking of a desired reference signal. It is shown that the proposed control strategy guarantees closed-loop stability with arbitrarily small steady-state errors. The  transient  performance  in the presence of non-zero initialization errors is quantified in terms of decreasing functions. Rigorous mathematical analysis and illustrative examples are provided to validate the theoretical claims.
\end{abstract}


\IEEEpeerreviewmaketitle

\section{Introduction} 
\label{sec:Introduction}

Adaptive control has been an active research topic in the past few decades, and it has been recognized as an effective approach to deal with systems that have uncertainties and disturbances \cite{Book_Sastry1989, Book_Narendra89, Book_KKK95, Book_Astrom1994, Book_L1}. Most of the success stories known to date have used state feedback approaches 
\cite{Calise_JGCD2000, Brinker_JGCD2001, Wise_GNC2005, Irene_GNC2011, Ackerman_GNC2016, Lee_JGCD2017}.
However, such approaches require that the state of the system is measurable, which is not always possible in practice. For this reason, there has been a significant effort to develop output feedback extensions.

The literature concerned with adaptive output feedback control is mainly focused on SISO systems or MIMO systems with strict structural requirements. References \cite{Book_Ioan96,Tao_IJAS1988} extend the results for SISO SPR systems to square MIMO systems. A modified interactor is introduced in order to relax the SPR assumption, thus increasing the applicability of the result to square MIMO systems with high relative degree. Similarly, \cite{Costa_ACC2002} borrows concepts and tools from 
\cite{Kokotovic_CDC1992, Jankovic_TAC1997}
to address square MIMO systems with arbitrary relative degree. One salient drawback in \cite{Costa_ACC2002} is the complicated structure of the controller, which makes it difficult to implement, especially as the relative degree increases. Nevertheless, the scope of these approaches is limited to square systems.

When dealing with \emph{non-square} MIMO systems, one common approach is to employ solutions for square systems in combination with squaring (-down or -up) methods.
Squaring-down methods can be applied to overactuated systems \cite{Sannuti_IJC1987} by reducing the excessive number of inputs. However, when dealing with underactuated systems, these methods discard some available measurements, thus limiting the use of output information. The disadvantage of squaring-down methods becomes even more evident when the system under consideration becomes non-minimum phase after squaring-down (e.g. missiles, inverted pendulums, etc.).

Recent work on adaptive output feedback control of under-actuated systems can be found in
\cite{Lavretsky_TAC2012, Gibson_TAC2015, Mizumoto_AUT2007, Lavretsky_AIAA2017}.
In particular, in \cite{Lavretsky_TAC2012, Gibson_TAC2015} solutions for square systems and their extensions to non-square systems are presented. These solutions are based on the use of the square-up method introduced in \cite{Misra_ACC1998}. These papers focus on systems in which the product between the input and output matrices is full rank. This assumption intrinsically implies that the system must have vector relative degree equal to $[1,\ldots,1]$, thus limiting the applicability of the approach. 
In \cite{Lavretsky_AIAA2017}, the authors augment the control law introduced in \cite{Lavretsky_TAC2012, Gibson_TAC2015} with a first order filter, thus extending the results to underactuated systems with arbitrary relative degree.
However, this approach assumes that the reference dynamics have vector relative degree $[1,\ldots,1]$. Moreover, the solution considers ideal parameterization of uncertainties  by an unknown constant matrix and known regressor functions. Thus, the work in \cite{Lavretsky_AIAA2017} does not lend itself to more general classes of non-square systems with time-varying uncertainties and unknown regressor functions, commonly found in many real-world systems.
Finally, in \cite{Mizumoto_AUT2007} the authors tackle non-square MIMO systems by designing an adaptive controller with \emph{multi-rate inputs}. Nevertheless, the approach requires the lifted system to be ASPR, and thus may not be applicable to systems with arbitrary relative degree.

In this paper, we propose an output feedback adaptive controller that deals with a general class of underactuated systems with arbitrary relative degree and with matched uncertainties.
The main contributions of this paper are: $(1)$ the controller handles underactuated MIMO systems with arbitrary relative degree and with time-varying uncertainties; $(2)$ uncertainties are not necessarily parameterized by known regressor functions, which broadens the applicability of the solution when compared to existing results; $(3)$ the approach is based on the right interactor matrix and a suitably defined state decomposition, providing semi-global stabilization for uncertain systems; $(4)$ the solution exhibits guaranteed performance during the transient and steady state under mild assumptions on the uncertainties and unknown initialization error.

The approach  is based on $\Lone$ adaptive control theory, which introduces a filtering structure providing a trade-off between robustness and performance. With this architecture, the filtering structure decouples the estimation loop from the control loop, thus allowing high-adaptation gains. While $\Lone$ adaptive state-feedback controllers (e.g. \cite{Cao_TAC_08Feb, Xargay_ACC2010}) have been successfully employed in real applications \cite{Jiang_JGCD2008, Griffin_GNC2010, Isaac_JGCD10, L1_Safety11,  Bichlmeier_GNC2013, Lee_JGCD2017, Ackerman_JGCD2017}, the literature directly concerned with output-feedback problems is less extensive. $\Lone$ output-feedback solutions for Single-Input Single-Output (SISO) systems can be found in \cite{L1_NonSPR09,Evgeny_ACC2011,Lee_ACC2014}, and can be easily extended to square MIMO systems \cite{Mahdianfar_JPC2016}.

An adaptive control solution for underactuated MIMO systems is presented in \cite{Lee_TAE2018}, where a suitably defined state decomposition is introduced, which enables standard $\Lone$ adaptive output feedback controllers to tackle underactuated systems. Nevertheless, the approach is limited to systems with relative degree one. The present article builds on and extends the work reported in \cite{Lee_TAE2018} to a more general class of systems with arbitrary relative degree by introducing modified $\Lone$ adaptive control laws based on the right interactor matrix.

This paper is organized as follows: in Section \ref{sec:MathematicalPreliminaries} we introduce mathematical results used in the paper; in Section \ref{sec:ProblemFormulation} a formal definition of the problem at hand is given; in Section \ref{sec:ControllerDesign} the main result of the paper is presented; Section \ref{sec:Stability Analysis} derives transient and steady-state performance of the system; in Section \ref{sec:Example} illustrative examples are provided to validate the theoretical findings; finally, the paper ends with concluding remarks in Section \ref{sec:Conclusions}.

\section{Mathematical preliminaries}\label{sec:MathematicalPreliminaries}
In this section we introduce few theoretical results that will be used in the paper. {Throughout the paper we use $\Vert \cdot \Vert$ to denote the vector or matrix $\infty$-norm. Given a signal $x(t)$, $\Linfnorm{x}$ and $\Linfnormtb{x}{a}{b}$ denote the $\Linf$ norms over $[0, \infty)$ and $[a,b]$, respectively. Finally, $\Linfnormt{x}{\tau}$ denotes the truncated norm $\Linfnormtb{x}{0}{\tau}$.
\begin{defn}
Let $M_0(s)$ be a $p \times m$ transfer matrix with $m \leq p$, and $S_{M_0}(s)$ be the Smith-McMillan form of $M_0(s)$. Suppose the normal rank of $M_0(s)$ is $m$. Let the polynomial $p_i(s)$ be the $i$-th diagonal element of $S_{M_0}(s)$, $i \in [1, \cdots, m]$. Then, the vector relative degree of $M_0(s)$ is defined as $[r_1, \cdots, r_m] \in \IR^m$, where $r_i$ is the relative degree of $p_i(s)$.
\end{defn}

\begin{defn}
Let $M_0(s)$ be be a $p \times m$ transfer matrix with $m \leq p$. Suppose $M_0(s)$ has the full normal column rank $m$. Then, $Z_0^{-1}(s)$ is called a right interactor of $M_0(s)$ if
$$
\lim\limits_{s \rightarrow \infty}( M_0(s) Z_0^{-1}(s) ) $$
is full rank.
\end{defn}

The following theorem is derived from \cite{Xin_LAA2002}.
\begin{thm}\label{thm:Interactor}
Let $M_0(s) = C \left( s\II_n - A \right)^{-1} B + D$, where $A \in \IR^{n \times n}$, $B \in \IR^{n \times m}$, $C \in \IR^{p \times n}$, and $D \in \IR^{p \times m}$. Assume that $(A,C)$ and $(A,B)$ are observable and controllable pairs, respectively. Suppose $M_0(s)$ has full normal column rank $m$ with $p\geq m$. Then, there exists a right interactor $Z_0^{-1}(s)$ such that
\begin{equation*}
Z_0(s) = C_z \left( s \II_{n_z} - {A}_z \right)^{-1} {B}_z + D_z,
\end{equation*}
with
\begin{equation}\label{eq:InvertibleInteractor}
\left|    \left[ \begin{array}{cc} -s\II_{n_z} +A_z & B_z \\ C_z & D_z \end{array} \right]  \right| \neq 0 , \quad \forall s \in \IC ,
\end{equation}
where $A_z \in \IR^{n_z \times n_z}$ is Hurwitz, $B_z \in \IR^{n_z \times m}$, $C_z \in \IR^{m \times n_z}$, and $D_z \in \IR^{m \times m}$.
Moreover, there exist $T_z \in \IR^{n \times n_z}$, $\bar{B} \in \IR^{n \times m}$, and $\bar{D} \in \IR^{p \times m}$ such that
\begin{equation}\label{eq:InteractorEquivalence}
\begin{split}
& \left[ \begin{array}{cc} -s\II_n +A & B \\ C & D \end{array} \right]
\left[ \begin{array}{cc} T_z & 0 \\ 0 & \II_m \end{array} \right] \\
 & \qquad =
\left[ \begin{array}{cc} T_z & \bar{B} \\ 0 & \bar{D} \end{array} \right]
\left[ \begin{array}{cc} -s\II_{n_z} +A_z & B_z \\ C_z & D_z \end{array} \right] , \quad \forall s \in \IC,
\end{split}
\end{equation}
where $T_z$, $\bar{D}$ are full (column) rank, and $(A,\bar{B})$ is a controllable pair  satisfying
\begin{equation}
\bar{M}_0(s) = C ( s\II_n - A )^{-1} \bar{B} + \bar{D} = M_0(s) Z_0^{-1}(s).
\end{equation}
\end{thm}
\begin{proof}
See \cite{Xin_LAA2002}.
\end{proof}

\begin{rem}
The right interactor $Z_0^{-1}(s)$ is not unique. In fact, the zeros of $Z_0^{-1}(s)$ are the eigenvalues of $A_z$, which can be arbitrarily chosen. As long as the intersection of $\text{eig}({A_z})$ and $\text{eig}({A})$ is an empty set, the controllability of $(A,\bar{B})$ is guaranteed. The reader is referred to \cite{Xin_LAA2002, Xin_IJC1998} for additional details on how to compute the interactor and the associated matrices ($T_z$, $\bar{B}$, and $\bar{D}$).
\end{rem}

Now, let $M(s)$ be the stable transfer matrix such as:
\begin{equation}\label{eq:GeneralSystem}
M(s) = C_m \left( s\II_n - A_m \right)^{-1} B_m ,
\end{equation}
where ${A_m \in \IR^{n \times n}}$, $B_m \in \IR^{n \times m}$, and $C_m \in \IR^{p \times n}$ are a minimal realization of $M(s)$ with $ m \leq p \leq n$.

\begin{cor}\label{cor:Interactor}
Consider the transfer matrix given in \eqref{eq:GeneralSystem}. Suppose $(C_mB_m)$ is rank deficient. Then,
there exist a stable transfer matrix $Z(s)$, and  matrices, $\bar{B} \in \IR^{n \times m}$, $T_z \in \IR^{n \times n_z}$ such that
\begin{equation}\label{eq:InteractorSS}
\begin{split}
Z(s) = & C_z ( s\II_{n_z} - A_z )^{-1} {B}_z + {D}_z, \\
\bar{M}(s) = & C_m ( s\II_n - A_m )^{-1} \bar{B} = M(s) Z^{-1}(s), \\
\end{split}
\end{equation}
and
\begin{equation}\label{eq:InteratorEquiv}
\begin{split}
A_m T_z =  & T_z {A}_z + \bar{B}C_z, \quad C_m A_m T_z = C_m \bar{B} C_z, \\
B_m = & T_z B_z + \bar{B} D_z,  \quad  C_mB_m =  C_m \bar{B} D_z,
\end{split}
\end{equation}
where $A_z \in \IR^{n_z \times n_z}$, $B_z \in \IR^{n_z \times m}$, $C_z \in \IR^{n_z \times m}$, and $D_z \in \IR^{n_z \times m}$ satisfy  \eqref{eq:InvertibleInteractor} and $T_z$ is of full column rank.
Moreover, the following hold:
\begin{itemize}
\item
$(A_m,\bar{B})$ is controllable, and $ (C_m \bar{B}) $ is full rank.
\item
If $M(s)$ has no unstable zeros, then $\bar{M}(s)$ does not possess unstable zeros.
\end{itemize}
\end{cor}
\begin{proof}
	The proof of Corollary \ref{cor:Interactor} is given in the Appendix.
\end{proof}

\begin{rem}
If $(C_mB_m)$ is full rank, then $Z(s) =\II_m$. Moreover, the rank condition on $(C_m B_m)$ is associated with the vector relative degree of MIMO systems. It can be easily shown that $(C_m B_m)$ is full rank, if and only if the vector relative degree is $\text{1}_m=[1,\ldots,1] \in \IR^m$. For systems with high relative degrees, $(C_mB_m)$ has rank deficiency.
\end{rem}

\begin{cor}\label{cor:InteactorAugmentedSystem}
Consider the state-space representation of the system \eqref{eq:GeneralSystem}:
\begin{equation*}
\begin{split}
\dot{x}(t) = & A_m x (t) + B_m u_x(t), \\
y(t) = & C_mx(t), \quad x(0)=x_0,
\end{split}
\end{equation*}
where $x(t) \in \IR^n$, $u_x(t) \in \IR^m$, $y(t) \in \IR^p$ are the state, input, and output vectors, respectively; $x_0 \in \IR^n$ is an initial condition. Let $x_v(t) \in \IR^{n}$ and $x_z(t) \in \IR^{n_z}$ be the states of the following cascaded system:
\begin{align}\label{eq:CascadeSystem}
\dot{x}_z(t) =  & A_z x_z(t) + {B}_z u_x(t), \quad u_v(t) =  C_z x_z(t) + D_z u_x(t) , \nonumber \\
\dot{x}_v(t) = & A_m x_v(t) + \bar{B} u_v(t), \quad  y_v(t) = C_m x_v(t),  \nonumber \\
x_v(0) =& x_0, \quad x_z(0) = 0,
\end{align}
where $y_v(t) \in \IR^p$ is the output vector, and $A_z \in \IR^{n_z \times n_z} , B_z \in \IR^{n_z \times m} , C_z \in \IR^{m \times n_z}, D_z \in \IR^{m \times m}$,  $\bar{B}\in \IR^{n \times m}$ are defined in Corollary \ref{cor:Interactor}.
Then, for all $t \geq 0$
\begin{equation}\label{eq:EquivStateOutput}
x(t) = x_v(t) + T_z x_z(t), \quad y_v(t) = y(t),
\end{equation}
where $T_z \in \IR^{n \times n_z}$ is full column rank satisfying \eqref{eq:InteratorEquiv}.
\end{cor}
\begin{proof}
	The proof of Corollary \ref{cor:InteactorAugmentedSystem} is  in the Appendix.
\end{proof}

\begin{rem}\label{rem:RelationWithTz}
Corollary \ref {cor:InteactorAugmentedSystem} provides a relationship between the states of the original system and the states of its cascaded representation.
\end{rem}

\begin{lem}\label{lem:ProjectionOperator}
Consider the system $M(s)$ given in \eqref{eq:GeneralSystem}. Suppose ${M}(s)$ has no unstable transmission zeros. From Corollary \ref{cor:Interactor}, let $\bar{M}(s) = C_m (s\II_n - A_m)^{-1} \bar{B}$ satisfying \eqref{eq:InteractorSS}. Define $H = \bar{B} \left( C_m\bar{B} \right)^{\dag}$, and $A_H= (\II_n-HC)A_m$, where $(C_m\bar{B})^{\dag}$ is the generalized left inverse of $\left(C_m\bar{B}\right)$. Then, $(\II_n-HC_m )\bar{B} =0$, and there exists a gain $K_v \in \IR^{n \times p}$ such that $A_v = \left( A_H + K_vC_m \right)$ is Hurwitz.
\end{lem}
\begin{proof}
	The proof of Lemma \ref{lem:ProjectionOperator} is  in the Appendix.
\end{proof}

The following remark will be used later in the Lyapunov analysis of the proposed adaptive controller.
\begin{rem}
Let $v(t) = \left(\II_n - H C_m \right)x_v(t)$, where $x_v(t) \in \IR^{n}$ is the state of the cascaded system \eqref{eq:CascadeSystem}. Then, $x_v(t) = v(t) + Hy(t)$ gives a state decomposition, where $y(t) \in \IR^{p}$ is the output. Since $\left( \II_n-HC_m \right)\bar{B} = 0$, the dynamics of $v(t)$ are not affected by the matched uncertainties.
\end{rem}

\section{Problem formulation}
\label{sec:ProblemFormulation}
\subsection{Problem statement}\label{sec:ProblemStatement}
Consider the following MIMO system
\begin{equation}\label{eq:nonlinearmodel1}
\begin{array}{l}
\dot x(t) = A_mx(t) + B_m( \omega u(t) + f(x,t) ) \,,\\
y(t) = C_mx(t) \,, \quad x(0)=x_0 \,,
\end{array}
\end{equation}
where $x(t) \in \IR^{n}$, $u(t) \in \IR^{m}$, $y(t) \in \IR^{p}$ are state, input and measurable output vectors, respectively, with $p \geq m$, and $x_0 \in \IR^{n}$ is an initial value. Moreover, ${A_m \in \IR^{n \times n}}$ is a known Hurwitz matrix, $B_m \in \IR^{n \times m}$ and $C_m \in \IR^{p \times n}$ are known matrices. Let ${(A_m, \, B_m, \, C_m)}$ be the minimal realization of $M(s) = C_m \left(s\II_n - A_m \right)^{-1}B_m$, which describes the desired dynamics of the closed-loop system; suppose $M(s)$ has full column rank $m$. Finally, $\omega >0$ is an unknown constant input gain, and $f: \IR \times \IR^n \rightarrow \IR^m$ is an unknown function representing matched uncertainties.

\begin{asm}\label{as:DesiredSystem}
$M(s)$ does not have unstable transmission zeros.
\end{asm}

\begin{asm}\label{as:UncertaintiesOfOmega}
The unknown constant input gain satisfies $\omega \in \Cset_\omega$, where $\Cset_\omega = [\omega_l, \omega_u]$ is a known compact set with $0< \omega_l < \omega_u$.
\end{asm}
\begin{asm}\label{as:UncertaintiesOfFunction}
    There exists $ b_0 >0$ such that
    \begin{equation*}
    \norm{ f(0, t) } < b_0 , \quad \forall t \geq 0,
    \end{equation*}
    where $b_0$ is a known constant.
    Moreover, for any $\delta > 0$ there exist $d_{\delta} > 0$, and $b_\delta  > 0$ such that
    \begin{equation*}
    {\left\| {\frac{{\partial f(x,t)}}{{\partial x}}} \right\|} \le {d_\delta } , \quad {\left\| {\frac{{\partial f(x,t)}}{{\partial t}}} \right\|} \le {b_\delta } ,
    \quad \forall \norm{x} < \delta,
    \end{equation*}
    where $d_{\delta}$ and $b_\delta$ are known constants.
\end{asm}

\begin{prob} \label{prob:probformulation}
Consider the system described by Equation \eqref{eq:nonlinearmodel1} satisfying Assumptions \ref{as:DesiredSystem}-\ref{as:UncertaintiesOfFunction}.
Design a feedback control law for $u(t)$ such that $y(t)$ tracks the desired response $y_m(t)$ both in transient and steady state, where $y_m(t)$ is the signal with the Laplace transform of $ y_m(s) = M(s)K_g r(s) $ with $K_g \in \IR^{m \times m_r}$ being a feed-forward gain, and $r(t) \in \IR^{m_r}$ being a reference signal.
\end{prob}

\subsection{Parametrization of uncertain function}\label{sec:UncertaintyParametrization}
\begin{lem}\label{lem:UncertainZeroParametrization}
Let $\tau >0$, and let $X(t) = [X_1^\top(t), X_2^\top(t)]^\top$ be a continuous and (piecewise) differentiable function, where $X_1(t) \in \IR^{n_1} $, $X_2(t) \in \IR^{n_2}$. 
Suppose that $\Linfnorm{\dot{X}_\tau}$ is finite.
Consider a nonlinear function $f(X,t)$ satisfying Assumption \ref{as:UncertaintiesOfFunction} and
\begin{equation*}
\norm{f(X,t)} < \bar{d}_X \norm{X_1(t)} + \bar{b}_X, \quad \Linfnormt{X}{\tau} \leq \rho_X, \quad  0 \leq t \leq \tau,
\end{equation*}
for some $\rho_X > 0$, $\bar{d}_X >0$ and $\bar{b}_X >0$,
Then, there exist continuous and (piecewise) differentiable $\theta(t)$ and $\sigma(t)$, such that
\begin{equation*}
f(X,t) = \theta(t) \norm{X_1(t)} + \sigma(t) \, , \quad  \forall t \in [0, \tau],
\end{equation*}
and
\begin{equation*}
\begin{split}
\norm{\theta(t)} \leq \bar{d}_X, \quad \norm{\dot{\theta}(t)} \leq \bar{l}_\theta, \\
\norm{\sigma(t)} \leq \bar{b}_X, \quad \norm{\dot{\sigma}(t)} \leq \bar{l}_\sigma ,
\end{split}
\end{equation*}
where $\bar{l}_\theta$, $\bar{l}_\sigma$ are computable finite bounds.
\end{lem}

\begin{proof}
See \cite[Lemma A.9, Lemma A.10]{Book_L1}.
\end{proof}
From Corollary \ref{cor:Interactor}, let $\{A_z, B_z, C_z, D_z\}$ be the set of system matrices of $Z(s)$ defined for $M(s)$, and $T_z \in \IR^{n \times n_z}$, $\bar{B} \in \IR^{n \times m}$ be matrices satisfying \eqref{eq:InteratorEquiv}.
Consider the following systems:
\begin{equation}\label{eq:ZsU}
\begin{split}
\dot{x}_u(t) = & A_z x_u(t) + B_z  u(t) ,  \\
u_z(t) = & C_z x_u(t) + D_z u(t) , \quad x_u(0) = 0 , \\
\end{split}
\end{equation}
and
\begin{equation}\label{eq:ZsF}
\begin{split}
\dot{x}_f(t) = & A_z x_f(t) + B_z f( T_g x_g + T_z x_f,t) , \, \, \, \, x_f(0) = 0, \\
\end{split}
\end{equation}
where $x_g(t) = [x_v^\top(t), x_u^\top(t)]^\top$, $T_g =  [\II_n, \omega T_z ]$, and $f(\cdot, t)$ satisfies Assumption \ref{as:UncertaintiesOfFunction}. The state $x_v(t) \in \IR^n$ is  governed by the following \emph{virtual} system:
\begin{equation}\label{eq:VirtualSystem}
\begin{split}
\dot{x}_v(t) = & A_m x_v(t) + \bar{B} ( \omega u_z (t) + \bar{f}(X, t) ),  \\
y_v(t) = & C_m x_v(t), \quad x_v(0) = x_0,
\end{split}
\end{equation}
where
\begin{equation}\label{eq:barUncertainFunc}
\bar{f}(X, t) = C_z x_f(t) + D_z f( T_g x_g + T_z x_f,t) ,
\end{equation}
with $X = [x_g^\top (t), x_f^\top (t)]^\top$.
By letting $x_z(t) = x_f(t) + \omega x_u(t)$, from Corollary \ref{cor:InteactorAugmentedSystem} and Equations \eqref{eq:ZsU}~-~\eqref{eq:barUncertainFunc} it follows that $x(t) = T_g x_g(t) + T_z x_f(t)$, and $y_v(t) = y(t)$ for any $t \geq 0$, where $x(t)$, $y(t)$ are  solutions of \eqref{eq:nonlinearmodel1}.

The following lemma gives a parameterization of the unknown function $\bar{f}(X, t)$.
\begin{lem}\label{lem:UncertaintyParameterization}
Consider the systems given by Equations \eqref{eq:ZsU}-\eqref{eq:barUncertainFunc}. Let $\tau > 0$, $\rho_x >0$, and $\rho_u >0$. Suppose $\Linfnormt{x}{\tau} \leq \rho_x$, and $\Linfnormt{u}{\tau} \leq \rho_u$, where $x(t) = T_g x_g (t) + T_z x_f(t)$. The function $\bar{f}(X, t) $ in \eqref{eq:barUncertainFunc} can be parameterized as follows
\begin{equation} \label{eq:NonlinearParam}
\bar{f}(X, t) = \theta(t) \norm{x_g(t)} + \sigma(t) , \quad  0 \leq t \leq \tau,
\end{equation}
and
\begin{equation}\label{eq:TheSigBound}
\begin{split}
\norm{\theta(t)} \leq \bar{d}_{\rho_x} \quad \norm{\dot{\theta}(t)} \leq \bar{l}_\theta , \\
\norm{\sigma(t)} \leq \bar{b}_{\rho_x} , \quad \norm{\dot{\sigma}(t)} \leq \bar{l}_\sigma ,
\end{split}
\end{equation}
where $\bar{l}_\theta$, $\bar{l}_\sigma$ are computable finite bounds, and $\bar{d}_{\rho_x}$, $\bar{b}_{\rho_x}$ are given by
\begin{equation}\label{eq:DBBound}
\begin{split}
\bar{d}_{\rho_x} = & \max_{\omega \in \Cset_\omega}{ ( \norm{C_z T_z^\dag T_g} + \norm{D_z} \norm{T_g} d_{\rho_x} ) }, \\
 \bar{b}_{\rho_x} = & \norm{C_z T_z^\dag} \rho_x + \norm{D_z} \Lonenorm{T(s)} d_{\rho_x}^2 \rho_x  \\
& + \norm{D_z}( \Lonenorm{T(s)} d_{\rho_x} + 1 )b_0 ,\\
\end{split}
\end{equation}
with $T(s) = T_z (s\II_{n_z} - A_z )^{-1}B_z$, and $b_0$, $d_{\rho_x}$ defined in Assumption \ref{as:UncertaintiesOfFunction}.
\end{lem}
\begin{proof}
	The proof of Lemma \ref{lem:UncertaintyParameterization} is given in the Appendix.
\end{proof}

\begin{rem}\label{rem:ConservativeBoundsOnTheta}
Let $\eta_t(t) = \bar{f}(X,t)$. The signal $\eta_t(t)$ can be viewed as the lumped matched uncertainty of the virtual system (see \eqref{eq:VirtualSystem}). From Lemma \ref{lem:UncertaintyParameterization} the unknown signal $\eta_t(t)$ is represented by time-varying uncertain signals $\theta(t)$ and $\sigma(t)$. The conservative bounds of $\theta(t)$ and $\sigma(t)$ are estimated by \eqref{eq:TheSigBound}, depending on the choice of the right interactor $Z(s)$. Notice that from \eqref{eq:ZsF} and \eqref{eq:barUncertainFunc} $\eta_t(t)$ can be seen as the uncertainty filtered by $Z(s)$, since $x(t) = T_g x_g(t) + T_z x_f(t)$ holds.
\end{rem}

\section{$\Lone$ adaptive controller design}\label{sec:ControllerDesign}

Let $\rho_0 >0$ be a given constant satisfying $\norm{x_0} \leq \rho_0$ with $x_0 \in \IR^n$ being an initial condition, and let $\bar{\gamma} > 0$  be an arbitrarily small constant.
For a given $\delta >0 $ let
\begin{equation}\label{eq:LipschitzFunction}
L_\delta =\frac{ \bar{\delta}(\delta)  }{  \delta  }   d_{ \bar{\delta}(\delta)} , \quad \bar{\delta}(\delta)  = \delta + \bar{\gamma} \, ,
\end{equation}
where $ d_{ \bar{\delta}(\delta)}$ is introduced in Assumption \ref{as:UncertaintiesOfFunction}.
Let $Z^{-1}(s)$  be a right interactor of $sM(s)$ such that
\begin{equation*} 
Z(s) = C_z(s\II_{n_z} - A_z )^{-1} B_z +D_z,
\end{equation*}
where $\{A_z \in \IR^{n_z \times n_z}, B_z \in \IR^{n_z \times m}, C_z \in \IR^{m \times n_z}\}$ is a minimal realization of $Z(s)$. Notice that the existence of $Z(s)$ is guaranteed by Corollary \ref{cor:Interactor}. Let $T_z \in \IR^{n \times n_z}$ and $\bar{B} \in \IR^{n \times m}$ be matrices that satisfy \eqref{eq:InteratorEquiv}.
Let
$K_v \in \IR^{n \times p}$ be a stabilizing gain so that
\begin{equation}\label{eq:DefAv}
A_v = A_H +K_vC_m
\end{equation}
is Hurwitz (from Lemma \ref{lem:ProjectionOperator} such $K_v$ exists), where
\begin{equation}\label{eq:DefAH}
A_H = (\II_n - H{C_m}){A_m}, \quad H = \bar{B} (C_m \bar{B} )^\dag ,
\end{equation}
with $(C_m \bar{B} )^\dag$ being the generalized inverse of $(C_m \bar{B} )$.
Let $P_y \in \IR^{p \times p}$ be a given positive definite matrix, and $P_v \in \IR^{n \times n}$ be the positive definite matrix, which solves
\begin{equation}\label{eq:LyapunovEquation}
A_v^\top P_v + P_v A_v = -Q
\end{equation}
for a positive definite $Q  \in \IR^{n \times n}$ with $\epsilon_q < \lambda_{\min}(Q)$.
Define
\begin{equation}\label{eq:DefKappa}
\begin{array}{c}
\kappa_m  = \sup\limits_{t \geq 0} \norm{ e^{A_m t} }, \\
\kappa_y  = \sqrt{ n \frac{  \lambda_{\max}( \bar{P}_v ) }{  \lambda_{\min}(P_v)  }}, \quad
\kappa_v   = \sqrt{ n \frac{   \lambda_{\max}( \bar{P}_v ) }{  \lambda_{\min}(P_y)  }},
\end{array}
\end{equation}
where $\bar{P}_v = (\II_n - HC_m)^\top P_v (\II_n - HC_m)$.
Let $D(s)$ be a $m \times m$ transfer matrix such that for all $\omega \in \Ccal_\omega$
\begin{equation*}\label{eq:DefCs}
C(s) =  \omega C_0(s)
\end{equation*}
is stable with $C(0) = \II_m$, and $C(s)Z^{-1}(s)$ is strictly proper, where
\begin{equation}\label{eq:DefC0}
C_0(s) = D(s) ( \II_m + \omega D(s) )^{-1} .
\end{equation}
Moreover, it is assumed that $D(s)$ ensures that there exists $\rho_r > 0$ such that
\begin{equation}\label{eq:L1Condition}
\Lonenorm{G(s)}  <  \frac{\rho_r - \rho_{ext} - {\rho}_{int} }{L_{{\rho}_{r}}  \rho_r } , \quad \omega \in \Cset_\omega,
\end{equation}
where
\begin{equation}\label{eq:RhoEZKappa}
\begin{split}
\rho_{ext}  = &\Lonenorm{H_r(s)} \Linfnorm{r} + \Lonenorm{G(s)} b_0 ,\\
\rho_{int}  = & ( \kappa_m + \kappa_x ) \rho_0, \\
\kappa_x   = &  \Lonenorm{H_1(s)}\kappa_y + \Lonenorm{H_2(s)}\kappa_v ,
\end{split}
\end{equation}
with $\kappa_m$, $\kappa_y$, and $\kappa_v$ being given in \eqref{eq:DefKappa}. Moreover,
\begin{align}\label{eq:DefHsfun}
H_r(s)  = & H_0(s) C(s) K_g, \quad
H_0(s)  =    (s \II_n - A_m )^{-1}B_m,   \nonumber \\
H_1(s)  = & \omega H_0(s) C_1(s) , \quad H_2(s) =  \omega H_0(s) C_2(s) , \\
G(s)    = & H_0(s) \left( \II_m -C(s) \right) , \nonumber
\end{align}
and
\begin{equation}\label{eq:DefCsfun1}
\begin{split}
C_1(s) =& (s+\alpha)C_0(s)Z^{-1}(s)(C_m \bar{B})^\dag , \\
C_2(s) =& C_0(s)Z^{-1}(s)(C_m \bar{B})^\dag C_m A_m , \\
\end{split}
\end{equation}
where $\alpha>0$ will be defined later.
Notice that $L_{\rho_r}$ satisfies \eqref{eq:LipschitzFunction} with $d_{\rho_x}$ and
\begin{equation}\label{eq:DefRhox}
\rho_x = \rho_r + \bar{\gamma}.
\end{equation}
Finally, let $\alpha>0$ be chosen to satisfy
\begin{equation}\label{eq:AlphaY}
\alpha_y = 2\alpha - \alpha_\phi>0, \quad \alpha_\phi= \frac{m \bar{d}_{\rho_x}^2 }{  \epsilon_q }{\left\|{\sqrt{P_y}{C_m}\bar{B}}  \right\|_2^2  },
\end{equation}
where $\bar{d}_{\rho_x}$ is given in \eqref{eq:DBBound}, and $\sqrt{P_y} \in \IR^{p \times p}$ is the upper triangular matrix satisfying the Cholesky decomposition; $ P_y = \sqrt{P_y}^\top \sqrt{P_y}$.

\begin{rem}
Clearly, for small $\bar{\gamma} >0$ we have $\rho_x \approx \rho_r$; $\rho_r$ is used to characterize the conservative bounds on the positively invariant set for the states of the closed-loop system.
\end{rem}

Consider the following control law
\begin{equation}\label{eq:ControlLaws}
u(s) = -D(s)Z^{-1}(s) ( \hat{\eta}_t(s) - r_z(s)  ) \, ,
\end{equation}
where $r_z(s) = Z(s)K_g r(s)$, and $\hat{\eta}_t(s)$ is the Laplace transform of
\begin{equation}\label{eq:DefHatEtaT}
\hat{\eta}_t(t) = \hat{\omega}(t) u_z(t) + \hat{\theta}(t) \norm{\hat{x}_g(t)} + \hat{\sigma}(t),
\end{equation}
and $\hat{\omega}(t)$, $\hat{\theta}(t)$, $\hat{\sigma}(t)$ are  the adaptive estimates, $u_z(t)$ is given in \eqref{eq:ZsU},  $x_g(t) = [\hat{x}_v^\top(t), x_u^\top(t)]^\top$, $x_u(t)$ is defined in \eqref{eq:ZsU}, and $\hat{x}_v(t) = \hat{v}(t) + Hy(t)$ with $\hat{v}(t)$ being given by the following predictor:
\begin{equation}\label{eq:StatePredictor}
\begin{split}
\dot{\hat{v}}(t) = & A_v \hat{x}_v(t) - K_v {y}(t)  - P_v^{-1} A_m^\top C_m^\top P_y \tilde{y}(t) , \\
\dot{\hat{y}}(t) = & -\alpha \tilde{y}(t) + C_m A_m \hat{x}_v(t) + C_m \bar{B} \hat{\eta}_t(t) , \\
\hat{v}(0) = & 0 , \quad \hat{y}(t) = y_0, \\
\end{split}
\end{equation}
where $y_0= C_m x_0$ is assumed to be known, $\tilde{y}(t) = \hat{y}(t) - y(t)$, and $A_v$ is given in \eqref{eq:DefAv}.
Consider the following adaptive laws:
\begin{equation}\label{eq:AdaptiveLaws}
\begin{split}
\dot{\hat{\omega}}(t) =&  \Gamma_\omega \Proj{ \hat{\omega}(t) }{ - u_z^\top(t) e_y(t)  },  \quad \hat{\omega}(0) =1, \\
\dot{\hat{\theta}}(t) =&  \Gamma_\theta \Proj{ \hat{\theta}(t) }{ - \norm{\hat{x}_g(t)} e_y(t) },  \quad \hat{\theta}(0) =0, \\
\dot{\hat{\sigma}}(t) =& \Gamma_\sigma \Proj{ \hat{\sigma}(t) }{ - e_y(t) },  \quad \hat{\sigma}(0) =0, \\
\end{split}
\end{equation}
where $\Gamma_\omega >0$, $\Gamma_\theta >0$, $\Gamma_\sigma >0$ are  adaptation gains, and $e_y(t) = \bar{B}^\top C_m^\top P_y \tilde{y}(t)$. $\Proj{\cdot}{\cdot}$ denotes the projection operator which is widely used in adaptive control; the operator provides smooth transition of the estimates on the apriori known boundary of  uncertainties (see \cite{Pomet_TAC1992}).

\section{Stability Analysis}\label{sec:Stability Analysis}
Consider the following closed-loop reference system
\begin{equation}\label{eq:L1RefSysSS}
\begin{split}
{{\dot x}_{ref}}(t) = & {A_m}{x_{ref}}(t) +  {B_m}\left( {\omega {u_{ref}}(t) + f(\xref,t) } \right) ,\\
{y_{ref}}(t) = & C_m^{}{x_{ref}}(t), \quad x_{ref}(0)=0,
\end{split}
\end{equation}
with
\begin{equation}\label{eq:L1RefControl}
{u_{ref}}(s) =   C_0(s)\left(  K_g r(s) - \eta _{ref}(s)    \right) ,
\end{equation}
where $x_{ref}(t) \in \IR^{n}$, $y_{ref}(t) \in \IR^{p}$ are the reference system states and outputs, respectively, $r(s)$ is the Laplace transform of the reference command $r(t) \in \IR^{m_r}$, $K_g \in \IR^{m \times m_r}$ is a feed-forward gain, and $C_0(s)$ is given in \eqref{eq:DefC0}. Moreover, $\eta_{ref}(s)$ is the Laplace transform of $f(\xref,t)$.

The closed-loop reference system in \eqref{eq:L1RefSysSS} and \eqref{eq:L1RefControl} 
represents \emph{the best achievable performance} of the $\Lone$ adaptive architecture \cite{Book_L1}. It is not implementable since it depends on the unknowns; it is used only for analysis purposes.

\begin{lem}\label{lm:BIBSBIBO}
    Consider the closed-loop reference system given in \eqref{eq:L1RefSysSS} and \eqref{eq:L1RefControl} and design constraints defined via \eqref{eq:LipschitzFunction}~-~\eqref{eq:DefRhox}. Then, for each $\omega \in \Ccal_\omega$ and $\tau >0$ the following bound holds
    	\begin{equation}\label{eq:RefBoundedX}
    		\begin{split}
    		\Linfnormt{\xref}{\tau} \leq \rho_{rx} ,
    		\end{split}
    	\end{equation}
    	where
    	\begin{equation}\label{eq:DefRhoXR}
    	\begin{split}
    	\rho_{rx}  = \rho_r - \frac{\rho_{int}}{ 1 - \Lonenorm{G(s)} L_{\rho_r}} > 0,
    	\end{split}
    	\end{equation}
    	with $\rho_{int}$, $G(s)$ given in \eqref{eq:RhoEZKappa} and \eqref{eq:DefHsfun}, respectively.
    	Moreover,
    	\begin{equation}\label{eq:RefBoundedU}
    	\Linfnormt{\uref}{\tau} \leq  \rho_{ru},
    	\end{equation}
    	where
    	\begin{equation}\label{eq:DefRhoUR}
    		\rho_{ru}  = \Lonenorm{C_0(s)} \left(  L_{\rho_r} \rho_{rx} +b_0\right) + \Lonenorm{C_0(s)K_g} \Linfnorm{r} ,
    	\end{equation}    	
    	with $C_0(s)$ defined in \eqref{eq:DefC0}.
\end{lem}
\begin{proof}
	Notice that from \eqref{eq:L1Condition} and \eqref{eq:DefRhoXR} one has
    \begin{equation}\label{eq:RhoXrIneq}
    \rho_{rx} > \rho_{ext} \geq  0 \, ,
    \end{equation}
    where $\rho_{ext}$ is defined in \eqref{eq:RhoEZKappa}.
    Substituting the control law given by Equation \eqref{eq:L1RefControl} into \eqref{eq:L1RefSysSS}, it follows that
    \begin{equation}\label{eq:TfReferenceSystem}
    \begin{split}
    {x_{ref}}(s) = & H_r(s)r(s)
       + G(s){\eta _{ref}}(s) \\
    {u_{ref}}(s) = & C_0(s)(K_g r(s) - \eta _{ref}(s) ) , \\
    {y_{ref}}(s) = & C_m x_{ref}(s),
    \end{split}
    \end{equation}
    where ${\eta _{ref}}(s)$ is the Laplace transform of $f(\xref,t)$, and $C_0(s)$, $\{H_r(s), G(s) \}$ are given in \eqref{eq:DefC0} and \eqref{eq:DefHsfun}, respectively. The resulting closed-loop reference system given by Equation \eqref{eq:TfReferenceSystem} is equivalent to the one in \cite[Chapter 2]{Book_L1}. Therefore, the rest of the proof follows from \cite[Chapter 2]{Book_L1}, and is omitted for the sake of brevity.

\end{proof}

Notice that the stability of the reference system can be guaranteed by designing a filter with high-bandwidth (see Equation \eqref{eq:L1Condition}). However, a  high bandwidth filter may lead to loss of robustness to time delay \cite{Book_L1}. The choice of a filter defines the trade-off between performance and robustness.

Differently from existing $\Lone$ adaptive state-feedback solutions, the present approach additionally requires a minimum order  filter (i.e., $C(s)Z^{-1}(s)$ is proper). Such condition is typical for output-feedback approaches. For example, the methods of \cite{L1_NonSPR09,Evgeny_ACC2011} require choosing a low-pass filter dependent upon the system's relative degree. Since the $\Lone$ reference system is identical to that of the existing $\Lone$ state-feedback, the problem of designing an appropriate filter can be tackled by existing optimal filter design techniques (e.g., see \cite{Kim_JOTA2013}).

\begin{rem}
Notice that the condition given in \eqref{eq:L1Condition} depends on the upper bound of the partial derivative of $f(x,t)$ which, in turn, depends on the unknown initial condition. Thus, the stability result in Lemma \ref{lm:BIBSBIBO} is semi-global. However, in the case where the uncertain function $f(x,t)$ has globally bounded partial derivatives (e.g. $d_{\delta} \equiv L$ for some constant $L >0$), the stability results become global (see the details in \cite[Chapter 3]{Book_L1}).
\end{rem}

Now, the closed-loop system stability is analyzed and the transient and steady-state performance bounds are derived. To demonstrate the stability of the closed-loop system with the proposed $\Lone$ control laws \eqref{eq:ControlLaws}-\eqref{eq:AdaptiveLaws}, we show that the difference between the closed-loop system and the ideal reference system is semi-globally bounded with arbitrarily small steady-state bounds. Moreover, we demonstrate that the transient performance errors due to non-zero initial conditions are bounded by strictly decreasing functions. Before stating the main results, we introduce a few variables of interest. Let
\begin{equation}\label{eq:GammaXU}
\begin{split}
\gamma_{u_0} =  & \Lonenorm{C_0(s)} L_{\rho_r} \gamma_{x_0}   \\
& + (\Lonenorm{C_1(s)} \kappa_y + \Lonenorm{C_2(s)} \kappa_v),  \\
\gamma_{x_0 }  = & \frac{ \kappa_x  + \kappa_m }{ 1 - \Lonenorm{G(s)} L_{\rho_r}} ,\\
\gamma_u =  &\Lonenorm{C_0(s)} L_{\rho_r} \gamma_x + \frac{ \Lonenorm{C_1(s)}}{ \sqrt{\lambda_{\min}(P_y)}} + \frac{ \Lonenorm{C_2(s)}}{ \sqrt{\lambda_{\min}(P_v)}},  \\
\gamma_x  = & \frac{ {\lambda_{\min}(P_y)}^{-\frac{1}{2}}\Lonenorm{H_1(s)}  + {\lambda_{\min}(P_v)}^{-\frac{1}{2}}\Lonenorm{H_2(s)} }{1 - \Lonenorm{G(s)} L_{\rho_r}} ,
\end{split}
\end{equation}
where $\set{\kappa_m, \kappa_y,\kappa_v}$, $\kappa_x$, $\set{H_1(s), H_2(s)}$, and $\set{C_1(s), C_2(s)}$ are given in \eqref{eq:DefKappa}, \eqref{eq:RhoEZKappa}, \eqref{eq:DefHsfun}, and \eqref{eq:DefCsfun1}, respectively.
Let ${\epsilon_\gamma} > 0$ be a small constant that verifies
\begin{equation}\label{eq:BarEpsilon}
\gamma_x {\epsilon_\gamma} < \bar{\gamma} , \quad \gamma_u {\epsilon_\gamma} < \bar{\gamma}, \quad \forall \omega \in \Cset_\omega.
\end{equation}
Finally, let $\rho_u$, $\rho_{dx}$, and $\rho_{du}$  be
\begin{equation}\label{eq:RhoDXU}
\begin{array}{c}
\rho_{u}  = \rho_{ru} + \rho_{du} , \\
\rho_{dx} = \gamma_{x_0} \rho_0 + \bar{\gamma}, \quad \rho_{du} = \gamma_{u_0} \rho_0 + \bar{\gamma},
\end{array}
\end{equation}
respectively, where $\rho_{ru}$ is defined in \eqref{eq:DefRhoUR}.

\begin{lem}\label{lem:EstimationErrorBounds}
    Consider the system given by Equation \eqref{eq:nonlinearmodel1} with control law defined in \eqref{eq:ControlLaws}-\eqref{eq:AdaptiveLaws}. Let $\tau > 0$ be a positive constant. If $\Linfnormt{x}{\tau} \leq \rho_x$ and $\Linfnormt{u}{\tau} \leq \rho_u$, then for all $t \in [0, \tau]$ the output-estimation error verifies
    \begin{equation}\label{eq:EstimationErrorBounds}
    \begin{split}
		\norm{\tilde{y}(t)} \leq & \kappa_y e^{-\frac{\lambda_1}{2} t }\norm{x_0} + \sqrt{ \frac{ \theta_1  }{ \lambda_{\min}(P_y)}  }  \frac{1}{\sqrt\Gamma} ,
    \end{split}
    \end{equation}
    where $\kappa_y$ is defined in \eqref{eq:DefKappa}, and
    \begin{equation}\label{eq:EstimationVariables}
    \begin{split}
    & \lambda_1 =  \min \left( \frac{\lambda_{\min}(Q_v) }{ \lambda_{\max}(P_v)}, \alpha_y  \right)  ,    \\
    & \theta_1 = \theta_0 + 4 m \frac{ \bar{d}_{\rho_x} \bar{l}_\theta + \bar{b}_{\rho_x} \bar{l}_\sigma }{\lambda_1} , \\
    & \theta_0  = 4 \left( \omega_u^2 + m \bar{d}_{\rho_x}^2 + m \bar{b}_{\rho_x}^2 \right) ,\\
    & \Gamma   = \min \left( \Gamma_\omega, \Gamma_\theta, \Gamma_\sigma\right) ,
    \end{split}
    \end{equation}
    with $Q_v = Q - \epsilon_q \II_n \succ 0$, $\alpha_y >0$  given in \eqref{eq:AlphaY}, and $\bar{d}_{\rho_x}$,  $\bar{l}_\theta$, $\bar{b}_{\rho_x}$, $ \bar{l}_\sigma$  satisfying \eqref{eq:TheSigBound}.
\end{lem}
\begin{proof}
	The proof of Lemma \ref{lem:EstimationErrorBounds} is given in the Appendix.
\end{proof}

Lemma \ref{lem:EstimationErrorBounds} states that the output estimation errors are exponentially convergent to a set, whose bound depends on both the upper bound of  $\theta_1$ and the adaptation gain $\Gamma$. Equation \eqref{eq:EstimationErrorBounds} implies that high values of adaptation gains achieve arbitrarily small estimation errors.

\begin{thm}\label{thm:ClosedLoopStability}
     Consider the closed-loop system with $\Lone$ adaptive output feedback controller defined via \eqref{eq:ControlLaws}-\eqref{eq:AdaptiveLaws}, subject to the design constraints in \eqref{eq:LipschitzFunction}-\eqref{eq:AlphaY}. Suppose the adaptation gains are chosen sufficiently high to satisfy
     \begin{equation}\label{eq:GammaCondition}
     \Gamma > \frac{\theta_1}{ {\epsilon_\gamma}^2  },
     \end{equation}
     where $\Gamma$, $\theta_1$ are defined in \eqref{eq:EstimationVariables}, and $\epsilon_\gamma$ satisfies \eqref{eq:BarEpsilon}. Then, the following upper bounds hold:
     \begin{equation}\label{eq:TransientBounds1}
     \Linfnorm{ \xref -x  } \leq \rho_{dx} , \quad \Linfnorm{\uref -u  } \leq \rho_{du} ,
     \end{equation}
     and
     \begin{equation}\label{eq:TransientBounds2}
     \begin{split}
     \Linfnorm{  \yref -y    } \leq  \norm{C_m} \rho_{dx} , \\
     \Linfnorm{ x } \leq  \rho_x , \quad \Linfnorm{u } \leq \rho_u .
     \end{split}
     \end{equation}
     Moreover, for each $\omega \in \Cset_\omega$ there exist positive constants of $\gamma_{dx}$ and $\gamma_{dy}$, and strictly decreasing functions of $\upsilon_{dx}(t)$ and $\upsilon_{dy}(t)$, such that for all $t \geq 0$
     \begin{equation}\label{eq:TransientBounds3}
     \begin{split}
     \norm{\xref(t) - x(t)} \leq & \upsilon_{dx}(t) \norm{x_0} + \frac{ \gamma_{dx} }{ \sqrt \Gamma } , \\
     \norm{\yref(t) - y(t)} \leq & \upsilon_{dy}(t) \norm{x_0} + \frac{ \gamma_{dy} }{ \sqrt \Gamma } .
     \end{split}
     \end{equation}
\end{thm}
\begin{proof}
	The proof of Theorem \ref{thm:ClosedLoopStability} is shown in the Appendix.
\end{proof}

Theorem \ref{thm:ClosedLoopStability} implies that tracking errors asymptotically converge to an invariant set that can be made sufficiently small via high adaptation gains. Notice that $\upsilon_{dx}(t)$ and $\upsilon_{dy}(t)$ in \eqref{eq:TransientBounds3} are independent of the adaptation gain, which is subject to the lower bound in \eqref{eq:GammaCondition}. Therefore, the transient performance due to non-zero initial conditions is quantified by strictly decreasing functions, and the steady-state errors can be arbitrarily reduced by increasing the adaptation gain.

\begin{rem}
	In the present section the closed-loop stability is analyzed under the assumption that $y_0$ is known. In the case when $y_0$ is not precisely measured due to sensor noise, one can easily derive similar stability results following the same proof, setting $\hat{y}(0)=0$ and $\bar{P}_v = P_v$. 	
\end{rem}

\section{Illustrative Examples}\label{sec:Example}
In this section, two examples are illustrated to validate our claims.
\subsection{Academic example}
Consider the nonlinear system \eqref{eq:nonlinearmodel1} with
\begin{align*}
 A_m = &\left[ \begin{array}{ccc} -2 & 0 & 1 \\ 1 & -5 & 2 \\ 1 & 0  & -5.5 \end{array} \right], \quad B_m =\left[ \begin{array}{ccc} 2 \\ 2.5  \\  -3  \end{array} \right], \\
 C_m = &\left[ \begin{array}{ccc} -5 & 10 & 5 \\ 1.25 & -1 & 0 \end{array} \right],
\end{align*}
where the unknown input gain is $\Omega=0.8 \in [0.7, 1.2]$, and the nonlinear uncertainty is set to be $f(x,t) = f_1(x,t)$ with
\begin{align*}
f_1(x,t) = & 0.11 \norm{x}_2^2 + 0.23 x_1 \tanh( 0.5 x_1) x_1  + 1.24 x_2 x_3  \\
& + 0.8 (1- e^{-0.7t}) + 2.0.
\end{align*}
For the $\Lone$ adaptive controller $D(s)$ and $Z(s)$ are selected as $D(s) = \frac{5}{s(s/11+1)}$, and $Z(s) = \frac{1}{s/4+1}$.

In simulations, we let $r(t) = 2 + 2 \sin3t$ and $x_0 = [-0.6, 0.6, -0.9]^\top$. Figure \ref{fig:Ex1_F1} illustrates the output trajectories and control inputs of the reference system and the closed-loop system for the adaptation gain $\Gamma=500$; the time-delay margin is numerically investigated and  is $0.45s$. Figure \ref{fig:Ex1_Gamma} demonstrates time histories of the tracking errors (i.e., $\norm{x_{ref}-x}$) and estimation errors (i.e., $\norm{\hat{y}-y}$) for different choices of adaptation gains. The  steady-state tracking errors are reduced with high  adaptation gain, and the transient errors are decreasing over time for the non-zero initial condition, as expected per analysis in Section \ref{sec:Stability Analysis}. Finally, we consider a different uncertainty by letting $f(x,t) = f_2(x,t)$ with
\begin{align*}
f_2(x,t) = & 0.15 \norm{x}_2^2 + 0.22 x_1 \tanh( 0.2 x_1) x_1  + 1.34 x_2 x_3  \\
& + 0.5 (1- e^{-1.1t}) + 1.8,
\end{align*}
and apply the same controller. The system outputs and the control signal are plotted in Figure \ref{fig:Ex1_F2}; asymptotic tracking of the reference outputs is achieved without any retuning of the controller.

\begin{figure}[h]
	\centering
	\includegraphics[width=0.80\linewidth]{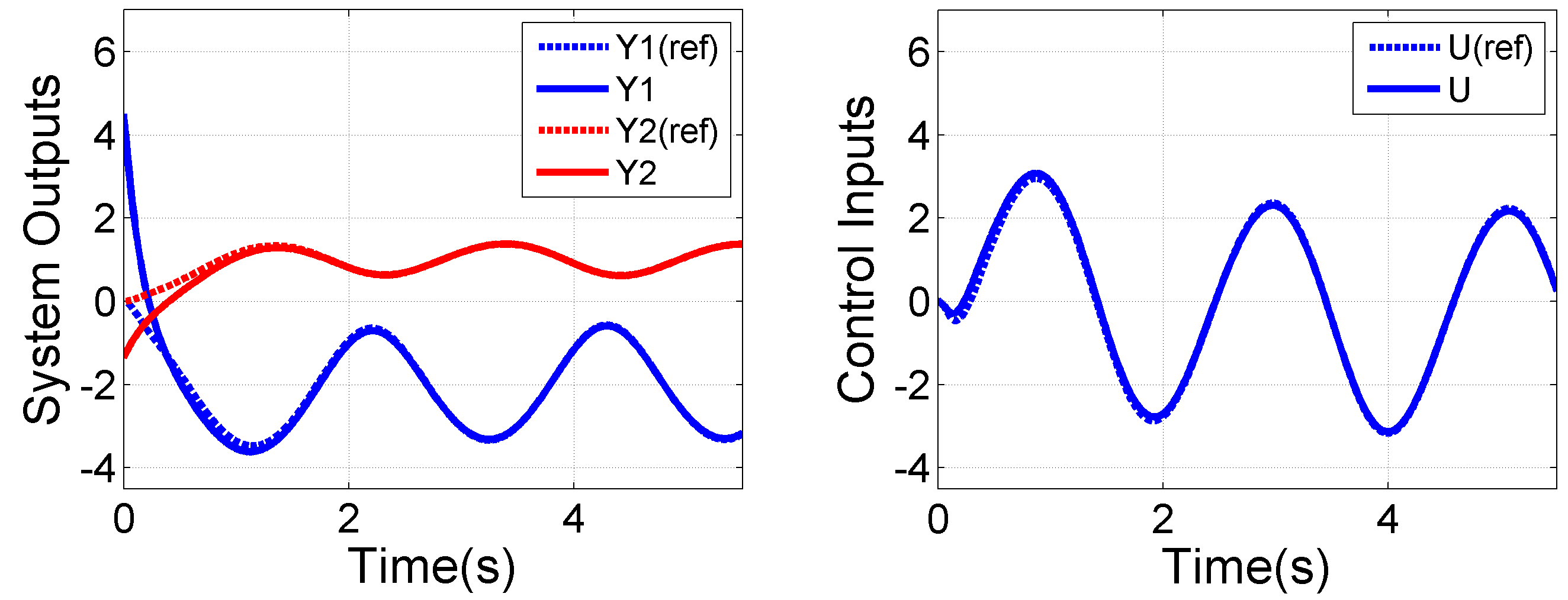}
	\caption{System responses with $\Gamma=500$ and $f_1(x,t)$}
	\label{fig:Ex1_F1}
\end{figure}
\begin{figure}[h]
	\centering
	\includegraphics[width=0.80\linewidth]{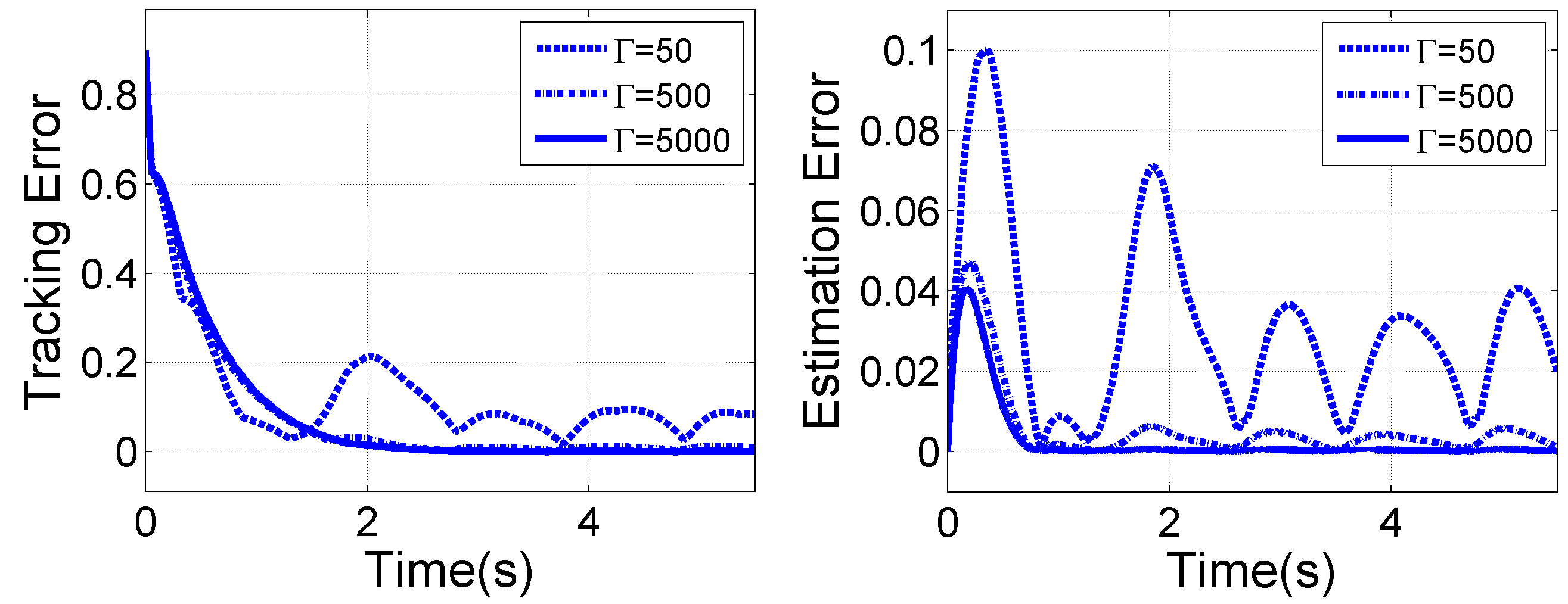}
	\caption{Tracking and estimation errors}
	\label{fig:Ex1_Gamma}
\end{figure}
\begin{figure}[h]
	\centering
	\includegraphics[width=0.80\linewidth]{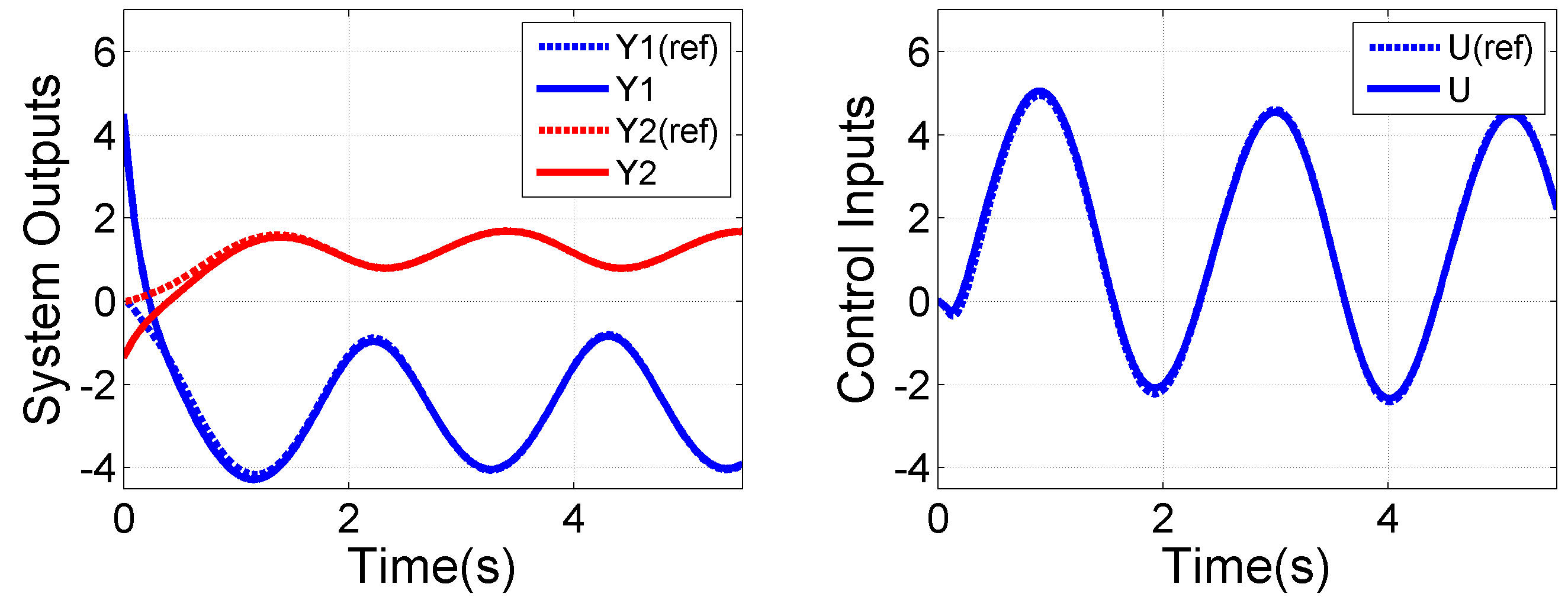}
	\caption{System responses with $\Gamma=500$ and $f_2(x,t)$}
	\label{fig:Ex1_F2}
\end{figure}

\subsection{Inverted pendulum on a cart}
In this section, we demonstrate the proposed method by designing the adaptive controller for an inverted pendulum on a cart. The control input is designed for the purpose of tracking a reference position, while maintaining the inverted pendulum balanced upright.
The nonlinear model of the inverted pendulum is given by
\begin{equation}\label{eq:nonlinear_PendulumOnCart}
 \begin{split}
\ddot{p}(t) - & \frac{\omega  u(t) - \nu \dot{p}(t) +F_{fric} (t) + d(t)}{M+m}  =   \\
& - \frac{m l \cos\theta(t) \ddot{\theta}(t) - m l  \sin \theta(t) \dot{\theta}^2(t)}{M+m} , \\
  m l  \cos \theta \ddot{p} (& t) - mg l  \sin \theta(t) + (I + ml ^2) \ddot{\theta}(t) = 0,
\end{split}
\end{equation}
where $p(t) \in \IR$, $\theta(t) \in \IR$ are the cart position and pendulum angle (measurable outputs), respectively;
$u(t)$ is the control input, $d(t)$ is the input disturbance, and $\omega >0$, $\nu>0$ are the motor constants.
$F_{fric}(t)$ is the nonlinear dynamic friction computed as \cite{Campbell_JDMC2008}:
\begin{equation}\label{eq:NonFriction}
\begin{split}
F_{fric}(t) = & - 73 \dot{p}(t) - 121 z(t) \left( 1 - 70 \frac{\norm{\dot{p}(t)}}{h(\dot{p}(t))} \right), \\
\dot{z}(t) = & \dot{p}(t) - 121 \frac{\norm{\dot{p}(t)}}{h(\dot{p}(t))} z(t),
\end{split}
\end{equation}
with $h(\dot{p}(t)) = - ( 0.04287 + 0.0432 e^{ -(\frac{\dot{p}(t)}{0.105} )^2 } ) (m+M)g$.
The nominal system parameters are selected as \cite{Campbell_JDMC2008}: $M_0 = 0.815$, $m_0 = 0.210$, $l_0 = 0.305$, $\omega_{0} = 1.719$, and $\nu_0 = 7.682$.
Moreover, it is assumed that the system has parameter variations from the nominal values, and therefore
\begin{equation}\label{eq:ParamVariation}
\begin{split}
M &=  1.2 M_0,  m = 0.8 m_0, \quad l = 1.2 l_0, \\
& \quad \omega = 1.2 \omega_{0}, \quad \nu  =1.5\nu_0 .
\end{split}
\end{equation}

For the purposes of comparison, we first consider a standard LQR controller for the system \eqref{eq:nonlinear_PendulumOnCart} \cite{Book_Panos1997}.  By letting $F_{fric}(t) \equiv 0$ and $d(t) \equiv 0$, the controller is obtained by linearizing the nonlinear model at $(p_e,\theta_e) =(0,0)$, together with $\cos\theta(t) \approx 1$. The LQR gain $K_{lqr}$ is given by
\begin{equation*}\label{eq:LQRgain}
K_{lqr} = [-7.0711 , -14.4505  ,-43.7667 ,  -7.6739].
\end{equation*}

For the $\Lone$ controller, the desired model is chosen  identical to the nominal (linearized) closed-loop system obtained by the LQR controller:
\begin{equation*}
\begin{split}
A_m = & \left[ \begin{array}{cccc}
0 & 1 & 0 & 0 \\ 14.62 & 20.64 & 88.23 & 15.87 \\ 0 & 0 & 0 & 1 \\ -44.26 & -62.47 & -237.34 & -48.04 \end{array} \right], \\
B_m = & \left[ \begin{array}{cccc} 0 & 2.07 & 0 & -6.26 \end{array} \right]^\top , \\
C_m = & \left[ \begin{array}{cccc} 1 & 0 & 0 & 0 \\ 0 & 0 & 1 & 0 \end{array} \right],
\end{split}
\end{equation*}
with the state vector $x(t)=[p(t),\dot{p}(t),\theta(t),\dot{\theta}(t)]^\top$, and the reference position command $r(t)$.
Since the desired model is obtained from the linearization, the uncertain function $f(x,t)$ in \eqref{eq:nonlinearmodel1} includes the linearization errors, parameter variations, non-linear friction $F_{fric}(t)$, and input disturbance $d(t)$. Notice that $(C_mB_m) = 0$. Therefore, we define the right interactor as $Z(s) = 0.47/(s+30)$, and choose $D(s) = \frac{30}{s (s/70+1) (s/100 +1)}$. The set of parameters of the $\Lone$ adaptive controller is given by $Q_v = 10\II_4$, $P_y = \II_2$, $\alpha = 25$, $\Gamma = 500$, and $K_g = -7.07$; the predictor gain $K_v$ is
given by
\begin{equation*}
K_v = \left[ \begin{array}{cccc}-4.51 & -22.087 & -1.56 & 36.98 \\ -1.56 & -22.91 & -2.87 & 40.55\end{array}  \right]^\top \, .
\end{equation*}

We present simulation results for two cases. The first case considers the nominal nonlinear dynamics, with $F_{fric}(t) \equiv 0$, $d(t) \equiv 0$, and zero initialization errors. Figure \ref{fig:Simuld1} illustrates and compares the system responses and control inputs
for the LQR controller and the $\Lone$ controller.
From the plots it can be noted that there is no significant difference in the performance of the solutions; this is not surprising, since the only uncertainties that affect the performance of the controllers are the linearization errors.
The second scenario considers the nonlinear system given by \eqref{eq:nonlinear_PendulumOnCart} with parametric variations given in \eqref{eq:ParamVariation}, with the nonlinear friction given by \eqref{eq:NonFriction}, input disturbance $d(t) = 3 sin(t)$, and non-zero initial conditions (the state is initialized as follows: $x_0 = [-0.5, -1, 0.1745, 0]^\top$). The results are illustrated in Figures \ref{fig:Simuld2} and \ref{fig:Simuld3}. As expected, the $\Lone$ controller ensures close tracking of the position, and boundedness of the angle within a neighborhood of zero, in spite of the uncertainties and non-zero initial errors, thus validating the theoretical claims.

\begin{figure}[h]
	\centering
	\includegraphics[width=0.8\linewidth]{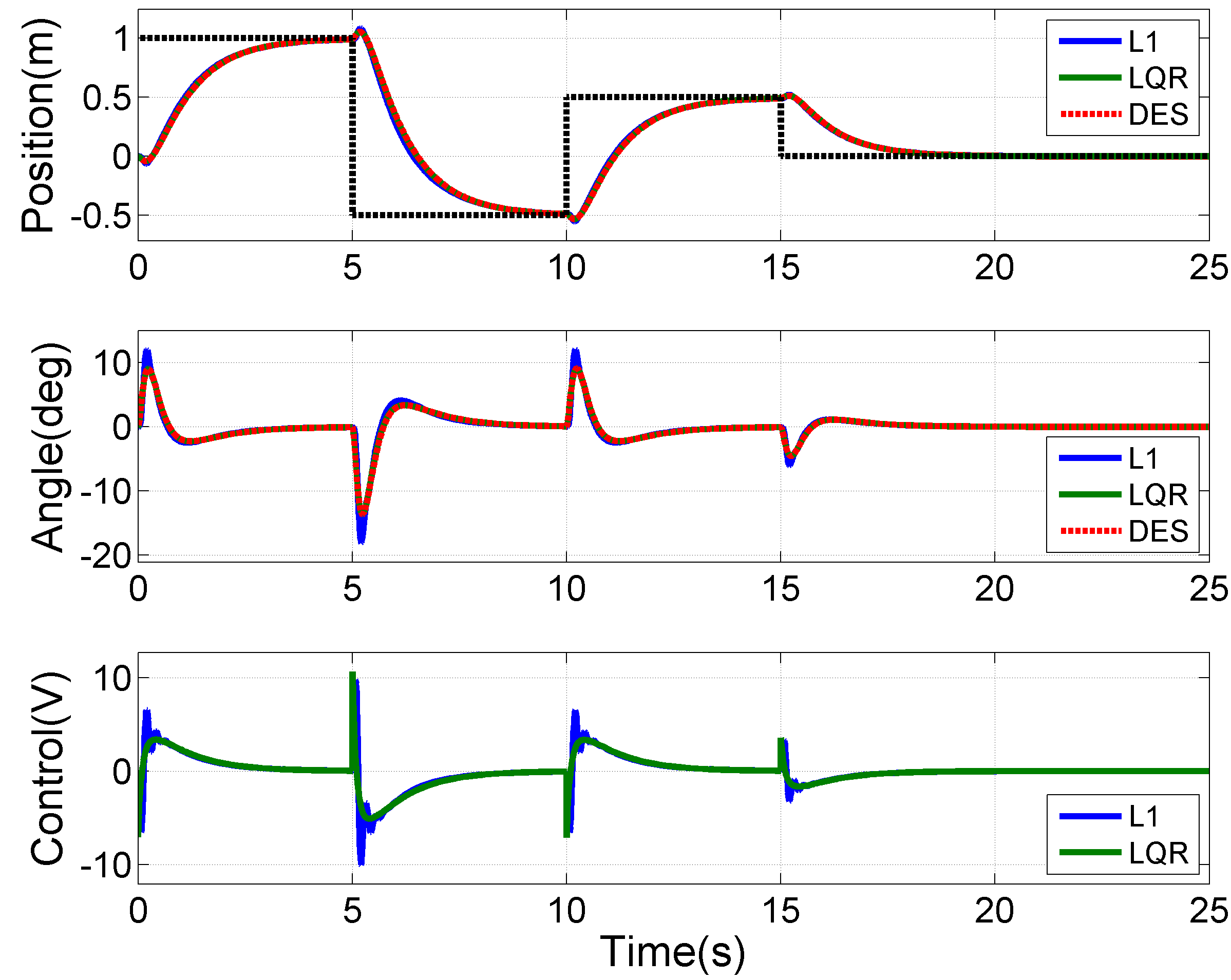}
	\caption{Inverted pendulum: position, angle, and control input for scenario 1.}
	\label{fig:Simuld1}
\end{figure}

\begin{figure}[h]
	\centering
	\includegraphics[width=0.8\linewidth]{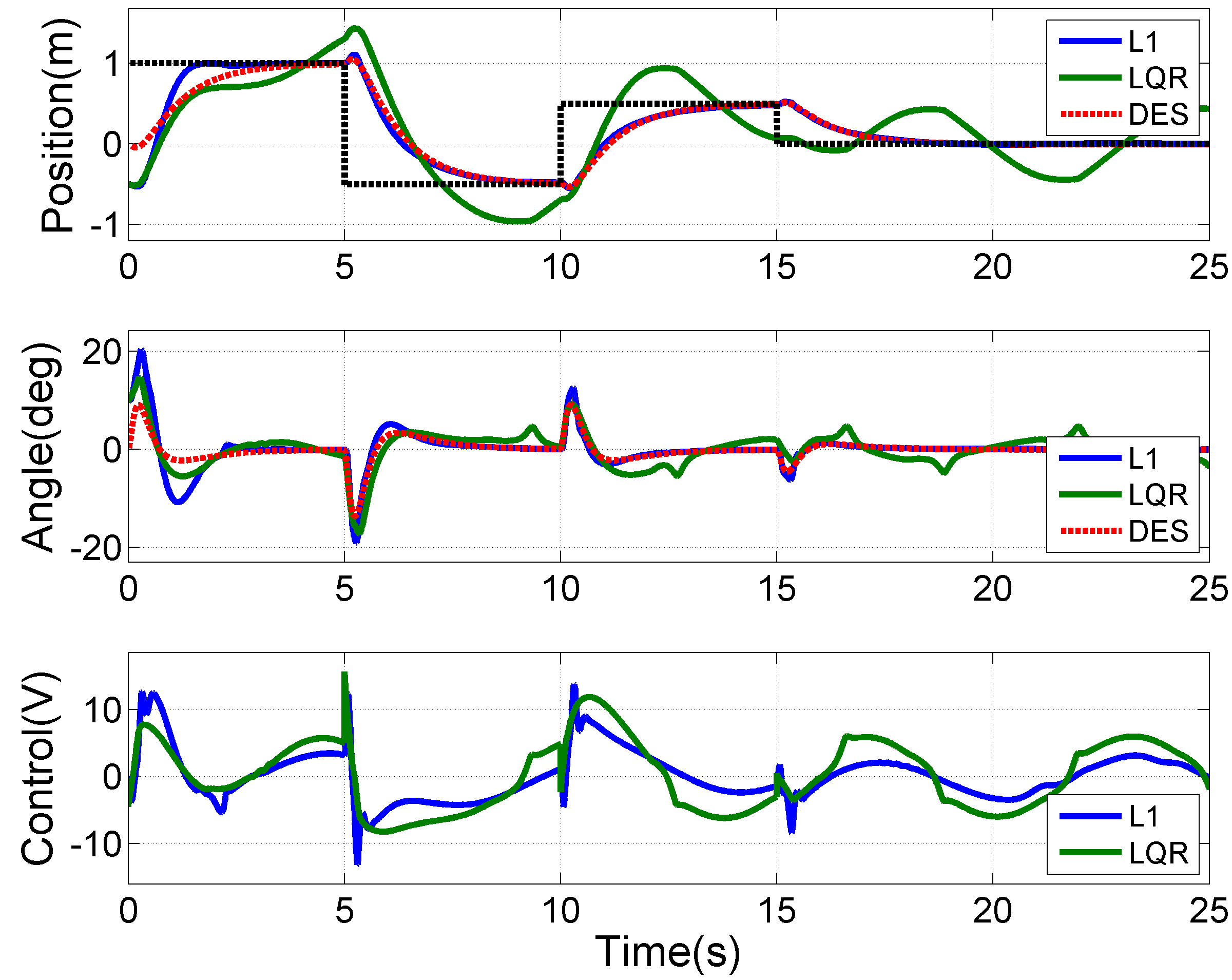}
	\caption{Inverted pendulum: position, angle, and control input for scenario 2.}
	\label{fig:Simuld2}
\end{figure}

\begin{figure}[h]
	\centering
	\includegraphics[width=0.85\linewidth]{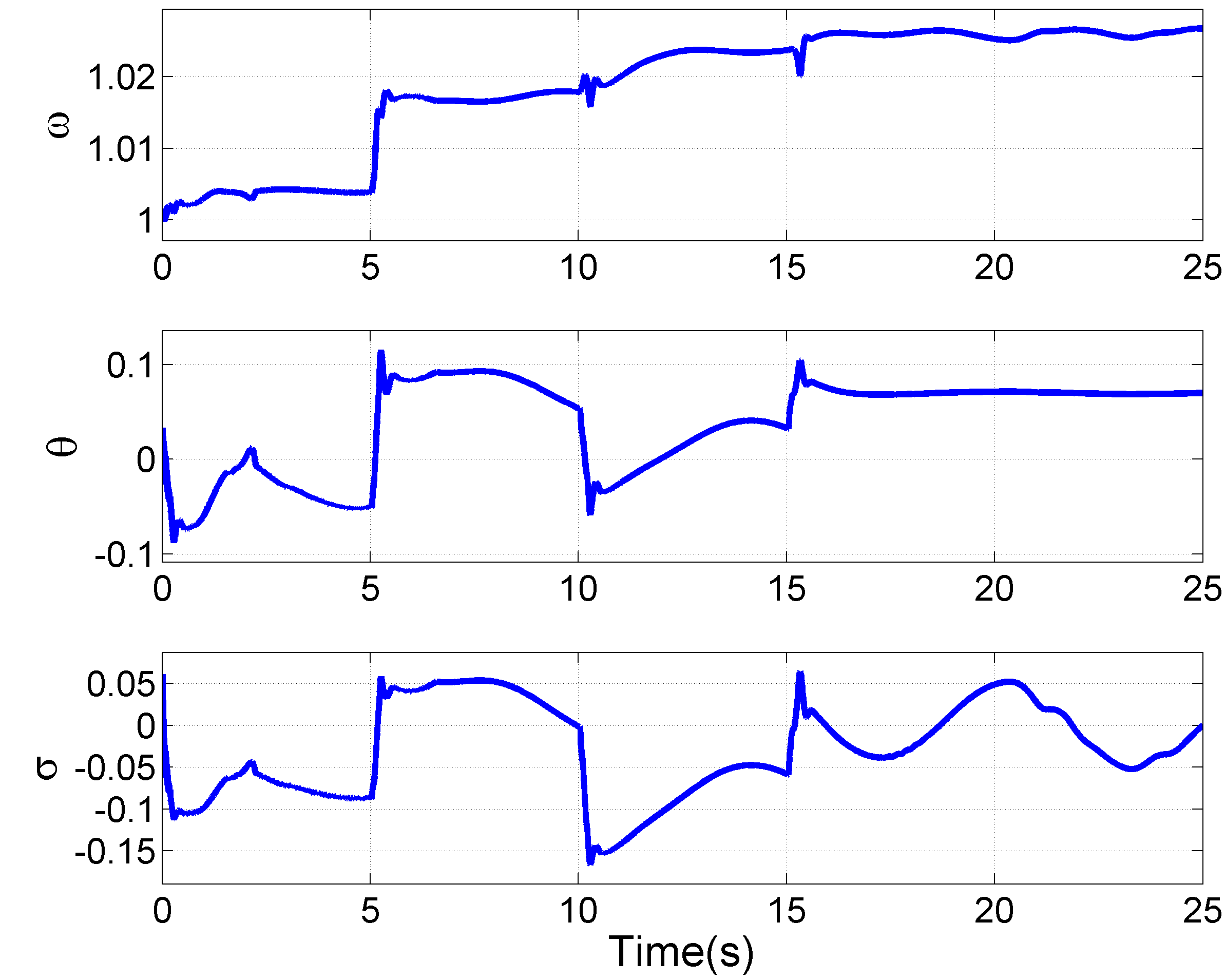}
	\caption{Inverted pendulum: adaptive gain dynamics for scenario 2.}
	\label{fig:Simuld3}
\end{figure}

\section{Conclusions} \label{sec:Conclusions}
This paper presents an $\Lone$ adaptive output feedback controller for non-square under-actuated MIMO systems with matched uncertainties. The controller design is based on the right interactor matrix, which is used to handle the non-square structure of the system through appropriate reparameterization of the system's equations. The control algorithm exhibits guaranteed performance in the transient and steady state under mild assumptions on the uncertainties and unknown initial error. Rigorous theoretical analysis and simulation results validate the performance of the proposed controller.

\bibliographystyle{ieeetr}
\bibliography{L1_TACBib}

\begin{thebibliography}{10}

\bibitem{Book_Sastry1989}
S.~Sastry and M.~Bodson, {\em Adaptive Control: Stability, Convergence and
  Robustness}.
\newblock {Advanced Reference}, Englewood Cliffs, NJ: {Prentice Hall}, 1989.

\bibitem{Book_Narendra89}
K.~S. Narendra and A.~M. Annaswamy, {\em Stable Adaptive Systems}.
\newblock {Information and System Sciences}, Englewood Cliffs, NJ: {Prentice
  Hall}, 1989.

\bibitem{Book_KKK95}
M.~Krsti{\'{c}}, I.~Kanellakopoulos, and P.~V. Kokotovi{\'{c}}, {\em Nonlinear
  and Adaptive Control Design}.
\newblock New York, NY: {John Wiley \& Sons}, 1995.

\bibitem{Book_Astrom1994}
K.~J. {\AA}str{\"{o}}m and B.~Wittenmark, {\em Adaptive control}.
\newblock Boston, MA: Addison-Wesley Longman Publishing Co., Inc., 1994.

\bibitem{Book_L1}
N.~Hovakimyan and C.~Cao, {\em $\mathcal{L}_1$ Adaptive Control Theory}.
\newblock Philadelphia, PA: {Society for Industrial and Applied Mathematics},
  2010.

\bibitem{Calise_JGCD2000}
A.~J. Calise, M.~Sharma, and J.~E. Corban, ``Adaptive autopilot design for
  guided munitions,'' {\em Journal of Guidance, Control, and Dynamics},
  vol.~23, pp.~837--843, September 2000.

\bibitem{Brinker_JGCD2001}
J.~S. Brinker and K.~A. Wise, ``Flight testing of reconfigurable control law on
  the x-36 tailless aircraft,'' {\em Journal of Guidance, Control and
  Dynamics}, vol.~24, pp.~903--909, September 2001.

\bibitem{Wise_GNC2005}
K.~Wise, E.~Lavretsky, J.~Zimmerman, J.~Francis, D.~Dixon, and B.~Whitehead,
  ``Adaptive flight control of a sensor guided munition,'' in {\em {AIAA}
  Guidance, Navigation, and Control Conference and Exhibit}, American Institute
  of Aeronautics and Astronautics ({AIAA}), August 2005.

\bibitem{Irene_GNC2011}
I.~M. Gregory, E.~Xargay, C.~Cao, and N.~Hovakimyan, ``Flight test of
  {$\mathcal{L}_1$}~adaptive control law: Offset landings and large flight
  envelope modeling work,'' in {\em AIAA Guidance, Navigation and Control
  Conference}, (Portland, OR), August 2011.
\newblock {AIAA}--2011--6608.

\bibitem{Ackerman_GNC2016}
K.~Ackerman, E.~Xargay, R.~Choe, N.~Hovakimyan, M.~C. Cotting, R.~B. Jeffrey,
  M.~P. Blackstun, T.~P. Fulkerson, T.~R. Lau, and S.~S. Stephens,
  ``$\mathcal{L}_1$ stability augmentation system for calspan's
  variable-stability learjet,'' in {\em AIAA Guidance, Navigation and Control
  Conference}, (San Diego, CA), January 2016.

\bibitem{Lee_JGCD2017}
H.~Lee, S.~Snyder, and N.~Hovakimyan, ``{$\mathcal{L}_1$} adaptive control
  within a flight envelope protection system,'' {\em Journal of Guidance,
  Control and Dynamics}, pp.~1--14, January 2017.

\bibitem{Book_Ioan96}
P.~A. Ioannou and J.~Sun, {\em Robust Adaptive Control}.
\newblock Upper Saddle River, NJ: {Prentice Hall}, 1996.

\bibitem{Tao_IJAS1988}
G.~Tao and P.~A. Ioannou, ``Robust model reference adaptive control for
  multivariable plants,'' {\em International Journal of Adaptive Control and
  Signal Processing}, vol.~2, pp.~217--248, September 1988.

\bibitem{Costa_ACC2002}
R.~Costa, L.~Hsu, A.~Imai, and G.~Tao, ``Adaptive backstepping control design
  for {MIMO} plants using factorization,'' in {\em American Control
  Conference}, (Anchorage, AK), {IEEE}, May 2002.

\bibitem{Kokotovic_CDC1992}
P.~Kokotovic, M.~Krstic, and I.~Kanellakopoulos, ``Backstepping to passivity:
  recursive design of adaptive systems,'' in {\em IEEE Conference on Decision
  and Control}, Institute of Electrical and Electronics Engineers ({IEEE}),
  December 1992.

\bibitem{Jankovic_TAC1997}
M.~Jankovic, ``Adaptive nonlinear output feedback tracking with a partial
  high-gain observer and backstepping,'' {\em IEEE Transactions on Automatic
  Control}, vol.~42, no.~1, pp.~106--113, 1997.

\bibitem{Sannuti_IJC1987}
P.~Sannuti and A.~Saberi, ``Special coordinate basis for multivariable linear
  systems{\textemdash}finite and infinite zero structure, squaring down and
  decoupling,'' {\em International Journal of Control}, vol.~45,
  pp.~1655--1704, May 1987.

\bibitem{Lavretsky_TAC2012}
E.~Lavretsky, ``Adaptive output feedback design using asymptotic properties of
  {LQG}/{LTR} controllers,'' {\em {IEEE} Transactions on Automatic Control},
  vol.~57, pp.~1587--1591, June 2012.

\bibitem{Gibson_TAC2015}
T.~E. Gibson, Z.~Qu, A.~M. Annaswamy, and E.~Lavretsky, ``Adaptive output
  feedback based on closed-loop reference models,'' {\em {IEEE} Transactions on
  Automatic Control}, vol.~60, no.~10, pp.~2728--2733, 2015.

\bibitem{Mizumoto_AUT2007}
I.~Mizumoto, T.~Chen, S.~Ohdaira, M.~Kumon, and Z.~Iwai, ``Adaptive output
  feedback control of general {MIMO} systems using multirate sampling and its
  application to a cart{\textendash}crane system,'' {\em Automatica}, vol.~43,
  pp.~2077--2085, December 2007.

\bibitem{Lavretsky_AIAA2017}
E.~Lavretsky, ``Robust and adaptive output feedback control for non-minimum
  phase systems with arbitrary relative degree,'' in {\em {AIAA} Guidance,
  Navigation, and Control Conference}, American Institute of Aeronautics and
  Astronautics, Jan 2017.

\bibitem{Misra_ACC1998}
P.~Misra, ``Numerical algorithms for squaring-up non-square systems, part ii:
  General case,'' in {\em 1998 American Control Conference}, (San Francisco,
  CA), June 1998.

\bibitem{Cao_TAC_08Feb}
C.~Cao and N.~Hovakimyan, ``Design and analysis of a novel
  $\mathcal{L}_1$~adaptive control architecture with guaranteed transient
  performance,'' {\em IEEE Transactions on Automatic Control}, vol.~53,
  pp.~586--591, March 2008.

\bibitem{Xargay_ACC2010}
E.~Xargay, N.~Hovakimyan, and C.~Cao, ``{$\mathcal{L}_1$} adaptive controller
  for multi-input multi-output systems in the presence of nonlinear unmatched
  uncertainties,'' in {\em American Control Conference}, Institute of
  Electrical and Electronics Engineers ({IEEE}), June 2010.

\bibitem{Jiang_JGCD2008}
J.~Wang, V.~V. Patel, C.~Cao, N.~Hovakimyan, and E.~Lavretsky, ``Novel
  $\mathcal{L}_1$ adaptive control methodology for aerial refueling with
  guaranteed transient performance,'' {\em Journal of Guidance, Control and
  Dynamics}, vol.~31, pp.~182--193, January--February 2008.

\bibitem{Griffin_GNC2010}
B.~Griffin, J.~Burken, and E.~Xargay, ``{$\mathcal{L}_1$} adaptive control
  augmentation system with application to the x-29 lateral/directional
  dynamics: A multi-input multi-output approach,'' in {\em AIAA Guidance,
  Navigation and Control Conference}, American Institute of Aeronautics and
  Astronautics, August 2010.

\bibitem{Isaac_JGCD10}
I.~Kaminer, A.~Pascoal, E.~Xargay, N.~Hovakimyan, C.~Cao, and V.~Dobrokhodov,
  ``Path following for unmanned aerial vehicles using $\mathcal{L}_1$~adaptive
  augmentation of commercial autopilots,'' {\em Journal of Guidance, Control
  and Dynamics}, vol.~33, pp.~550--564, March--April 2010.

\bibitem{L1_Safety11}
N.~Hovakimyan, C.~Cao, E.~Kharisov, E.~Xargay, and I.~M. Gregory,
  ``{$\mathcal{L}_1$} adaptive control for safety-critical systems,'' {\em IEEE
  Control Systems Magazine}, vol.~31, pp.~54--104, October 2011.

\bibitem{Bichlmeier_GNC2013}
M.~Bichlmeier, F.~Holzapfel, E.~Xargay, and N.~Hovakimyan, ``{$\mathcal{L}_1$}
  adaptive augmentation of a helicopter baseline controller,'' in {\em AIAA
  Guidance, Navigation and Control Conference}, American Institute of
  Aeronautics and Astronautics, August 2013.

\bibitem{Ackerman_JGCD2017}
K.~A. Ackerman, E.~Xargay, R.~Choe, N.~Hovakimyan, M.~C. Cotting, R.~B.
  Jeffrey, M.~P. Blackstun, T.~P. Fulkerson, T.~R. Lau, and S.~S. Stephens,
  ``Evaluation of an {$\mathcal{L}_1$} adaptive flight control law on calspan's
  variable-stability learjet,'' {\em Journal of Guidance, Control and
  Dynamics}, vol.~40, pp.~1051--1060, April 2017.

\bibitem{L1_NonSPR09}
C.~Cao and N.~Hovakimyan, ``$\mathcal{L}_1$ adaptive output-feedback controller
  for non-stricly-positive-real reference systems: Missile longitudinal
  autopilot design,'' {\em {AIAA Journal of Guidance, Control, and Dynamics}},
  vol.~32, pp.~717--726, May-June 2009.

\bibitem{Evgeny_ACC2011}
E.~Kharisov and N.~Hovakimyan, ``$\mathcal{L}_1$ adaptive output feedback
  controller for minimum phase systems,'' in {\em American Control Conference},
  (San Francisco, CA), June--July 2011.

\bibitem{Lee_ACC2014}
H.~Lee, V.~Cichella, and N.~Hovakimyan, ``$\mathcal{L}_1$ adaptive output
  feedback augmentation of model reference control,'' in {\em American Control
  Conference}, (Portland, OR), June 2014.

\bibitem{Mahdianfar_JPC2016}
H.~Mahdianfar, N.~Hovakimyan, A.~Pavlov, and O.~M. Aamo, ``{$\mathcal{L}_1$}
  adaptive output regulator design with application to managed pressure
  drilling,'' {\em Journal of Process Control}, vol.~42, pp.~1--13, June 2016.

\bibitem{Lee_TAE2018}
H.~Lee, S.~Snyder, and N.~Hovakimyan, ``{$\mathcal{L}_1$} adaptive output
  feedback augmentation for missile systems,'' {\em {IEEE} Transactions on
  Aerospace and Electronic Systems}, vol.~54, pp.~680--692, April 2018.

\bibitem{Xin_LAA2002}
X.~Xin and T.~Mita, ``A simple state-space design of an interactor for a
  non-square system via system matrix pencil approach,'' {\em Linear Algebra
  and its Applications}, vol.~351-352, pp.~809--823, August 2002.

\bibitem{Xin_IJC1998}
X.~Xin and T.~Mita, ``Inner-outer factorization for non-square proper functions
  with infinite and finite j omega -axis zeros,'' {\em International Journal of
  Control}, vol.~71, pp.~145--161, January 1998.

\bibitem{Pomet_TAC1992}
J.-B. Pomet and L.~Praly, ``Adaptive nonlinear regulation: Estimation from the
  {L}yapunov equation,'' {\em IEEE Transactions on Automatic Control}, vol.~37,
  pp.~729--740, June 1992.

\bibitem{Kim_JOTA2013}
K.-K.~K. Kim and N.~Hovakimyan, ``Multi-criteria optimization for filter design
  of $\mathcal{L}_1$ adaptive control,'' {\em Journal of Optimization Theory
  and Applications}, vol.~161, pp.~557--581, September 2013.

\bibitem{Campbell_JDMC2008}
S.~A. Campbell, S.~Crawford, and K.~Morris, ``Friction and the inverted
  pendulum stabilization problem,'' {\em Journal of Dynamic Systems,
  Measurement, and Control}, vol.~130, no.~5, p.~054502, 2008.

\bibitem{Book_Panos1997}
P.~J. Antsaklis and A.~N. Michel, {\em Linear systems}.
\newblock Englewood Cliffs, NJ: {McGraw-Hill}, 1997.

\bibitem{Davison_AUT1974}
E.~Davison and S.~Wang, ``Properties and calculation of transmission zeros of
  linear multivariable systems,'' {\em Automatica}, vol.~10, pp.~643--658,
  December 1974.

\end{thebibliography}

\appendix

\begin{proof}[\textbf{Proof of Corollary \ref{cor:Interactor}}]
	Notice that $(C_mB_m) = \lim\limits_{ s \rightarrow \infty} ( sM(s) )$. Let $M_0(s) = sM(s)$, and $M_0(s) = C ( s\II_n - A )^{-1} B + D$ with $A = A_m$, $B = B_m$, $C = C_mA_m$, and $D = C_mB_m$. Since $\{A_m, B_m, C_m\}$ is controllable-observable, and $A_m$ is Hurwitz, the triple $\{ A,B,C\}$ is also controllable-observable. Therefore, from Theorem \ref{thm:Interactor} it follows that there exists a right interactor $Z^{-1}(s)$ which satisfies \eqref{eq:InteractorEquivalence} with $T_z \in \IR^{n \times n_z}$, $\bar{B} \in \IR^{n \times m}$, and $\bar{D} \in \IR^{p \times m}$; $(A_m, \bar{B})$ is controllable. Since Equation \eqref{eq:InteractorEquivalence} holds, one has
	\begin{equation*} 
	\begin{split}
	(\bar{D} - C_m \bar{B} ) C_z  A_z^{-1} =  C_m T_z,  \quad
	(\bar{D} - C_m \bar{B} ) D_z =  C_m T_z B_z,
	\end{split}
	\end{equation*}
	which further leads to
	\begin{equation}\label{eq:CBprocedure2}
	(\bar{D} - C_m \bar{B} )( D_z -C_z  A_z^{-1}B_z ) = 0.
	\end{equation}
	Notice that both $(D_z - C_z A_z^{-1} B_z)$ and $\bar{D}$ are full rank (see Theorem \ref{thm:Interactor}).
	From \eqref{eq:CBprocedure2} it follows that $ \bar{D} = C_m \bar{B} $ holds. Therefore, $(C_m \bar{B})$ is full rank, and Equation \eqref{eq:InteratorEquiv}  follows from \eqref{eq:InteractorEquivalence}.
	
	Finally, suppose that $M(s)$ has no unstable transmission zeros. Notice that pole-zero cancellations in $\bar{M}(s)Z(s)$ only happen in $\IC^{-}$, since $A_z$ is Hurwitz, and therefore, $\bar{M}(s)$ cannot have any unstable transmission zeros.
	This completes the proof.
\end{proof}

\begin{proof}[\textbf{Proof of Corollary \ref{cor:InteactorAugmentedSystem}}]
	Notice that Equation \eqref{eq:CascadeSystem} can be rewritten as
	\begin{align}\label{eq:CascadeSS}
	\left[ \begin{array}{c}  \dot{x}_v(t) \\  \dot{x}_z(t) \end{array}  \right]
	& =
	\left[ \begin{array}{cc}  A_m & \bar{B} C_z \\ 0 & A_z \end{array} \right]
	\left[\begin{array}{c}  {x}_v(t) \\  {x}_z(t) \end{array}  \right]
	+
	\left[\begin{array}{c}  \bar{B} D_z \\  B_z  \end{array}  \right]  u_x(t) , \nonumber \\
	y_v(t) & = \left[  \begin{array}{cc} C_m & 0 \end{array} \right] \left[\begin{array}{c}  {x}_v(t) \\  {x}_z(t) \end{array}  \right].
	\end{align}
	Now, let $[ x_t^\top(t), x_z^\top(t)]^\top = T_t [ x_v^\top(t), x_z^\top(t)]^\top $ with $T_t =  \left[ \begin{array}{cc} \II_n & T_z \\ 0 & \II_{n_z} \end{array} \right]$.
	By applying a similarity transform with $T_t$, from \eqref{eq:CascadeSS} it follows, together with \eqref{eq:InteratorEquiv}, that
	\begin{align*} 
	\dot{x}_z(t) = & A_z x_t(t) + B_z u_x(t),  \\
	\dot{x}_t(t) = & A_m x_t(t) + B_m u_x(t),\quad y_v(t) =  C_m x_t(t),
	\end{align*}
	with $x_t(0) = x_0$, and $x_z(0) =0$. This implies that $x(t) = x_t(t)$ for all $t \geq 0$. Therefore, Equation \eqref{eq:EquivStateOutput} holds, which completes the proof.
\end{proof}

\begin{proof}[\textbf{Proof of Lemma \ref{lem:ProjectionOperator}}]
	Since $M(s)$ does not have any transmission zeros by hypothesis, from Corollary \ref{cor:Interactor} it follows that $\bar{M}(s)$ has no unstable transmission zeros, and $(C_m \bar{B})$ is full rank. Notice that $(C_m\bar{B})^\dag (C_m\bar{B}) = \II_m$, and therefore $H=B_m(C_m\bar{B})^\dag$ satisfies $(\II_n-HC_m )\bar{B} =0$.
	
	Next we show that $\left(A_H,C_m \right)$ is a detectable pair, which  ensures the existence of a stabilizing gain $K_v$ such that $A_v$ is Hurwitz. Suppose that ${z _i} \in \IC$  is an unobservable mode of $\left(A_H,C_m \right)$. By Popov-Belevitch-Hautus observability test\cite[Chapter 3]{Book_Panos1997}, there exists a non-zero vector  ${\xi_i} \in {R^n}$ such that
	\begin{equation}\label{eq:UnobservableVector1}
	A_H{\xi_i} = {z_i}{\xi_i}, \quad
	{C_m}{\xi_i} = 0,
	\end{equation}
	which yields
	\begin{equation*}\label{eq:UnobservableVector2}
	\left( {{z_i}\II_n - {A_m}} \right){\xi _i} + H{C_m}{A_m}{\xi_i} = 0.
	\end{equation*}
	Let $ {\varsigma_i} = {({C_m}\bar{B})^\dag }{C_m}{A_m}{\xi_i} \in \IR^{m} $. Then, it follows that
	\begin{equation}\label{eq:UnobservableVector3}
	{\bar{B}}{\varsigma_i} = {\bar{B}}{({C_m}{\bar{B}})^\dag }{C_m}{A_m}{\xi_i} = H{C_m}{A_m}{\xi_i}.
	\end{equation}
	Combining \eqref{eq:UnobservableVector1}~-~\eqref{eq:UnobservableVector3}, one has
	\begin{equation*} 
    R(z_i) \left[ {\begin{array}{*{20}{c}}
	{{\xi_i}}\\
	{{\varsigma_i}}
	\end{array}} \right] = \left[ {\begin{array}{*{20}{c}}
	0\\
	0
	\end{array}} \right] \,, \quad 	
	R(z_i) =  \left[ {\begin{array}{*{20}{c}}
		{{z _i}\II_n - {A_m}}&{\bar{B}}\\
		{ - {C_m}}&0
		\end{array}} \right],
	\end{equation*}
where $[\xi_i^\top,\varsigma_i^\top] \neq 0$. Since $\set{A_m, \bar{B}, C_m}$ is a minimal realization, from $\text{rank}(R(z_i)) < n + m$ it follows that $z_i$ must be the transmission zero of $M(s)$ \cite{Davison_AUT1974}. Finally, since $\bar{M}(s)$ does not possess unstable transmission zeros, $z_i \in \IC^{-}$ holds. Therefore, $\left( A_H,C_m \right)$ is detectable, which completes the proof.
\end{proof}

\begin{proof}[\textbf{Proof of Lemma \ref{lem:UncertaintyParameterization}}]
	Since $\Linfnormt{x}{\tau} \leq \rho_x$, from \eqref{eq:ZsF} it follows that
	\begin{equation}\label{eq:BoundTzxz}
	\Linfnormtb{T_z x_z}{0}{\tau} \leq  \Lonenorm{T(s)} d_{\rho_x} \rho_x + \Lonenorm{T(s)} b_0,
	\end{equation}
	which yields $\Linfnormt{X}{\tau} \leq \rho_X$ for some $\rho_X >0$, where $d_{\rho_x}$, $b_0$ are given in Assumption \ref{as:UncertaintiesOfFunction}, and
	\begin{equation*} 
	T(s) = T_z (s \II_{n_z} - A_z)^{-1} B_z.
	\end{equation*}
	Moreover, notice that $x_z(t) = T_z^\dag T_g x_g(t) + T_z^\dag x(t)$, where $T_z^\dag$ is the generalized inverse of $T_z$. From \eqref{eq:barUncertainFunc}, one has
	\begin{equation}\label{eq:BarUncertainFcn}
	\begin{split}
	\bar{f}(X, t) = & C_z T_z^\dag T_g x_g(t) + C_z T_z^\dag x(t) \\
	& + D_z f( T_g x_g + T_z x_f,t) .
	\end{split}
	\end{equation}
	Notice that using Assumption \ref{as:UncertaintiesOfFunction} on $f(x,t)$ it is easy to show that the partial derivatives of  $\bar{f}(X,t)$ are (semi-globally) bounded. Moreover, by using the fact that $\Linfnormtb{T_g x_g + T_z x_f}{0}{\tau} \leq \rho_x$, from \eqref{eq:BoundTzxz} and \eqref{eq:BarUncertainFcn} it follows that
	\begin{equation*}
	\begin{split}
	\norm{ \bar{f}(X, t) }& < \bar{d}_{\rho_x} \norm{x_g(t) } + \bar{b}_{\rho_x} , \quad 0 \leq t \leq \tau,
	\end{split}
	\end{equation*}
	where $\bar{d}_{\rho_x}$, $ \bar{b}_{\rho_x} $ are given in \eqref{eq:DBBound}.
	Notice that $\Linfnormt{x}{\tau} \leq \rho_x$ and $\Linfnormt{u}{\tau} \leq \rho_u$, along with \eqref{eq:ZsU}~-~ \eqref{eq:VirtualSystem}, imply that $\norm{\dot{X}(t)}$ is finite for all $0 \leq t \leq \tau$. Therefore, from Lemma \ref{lem:UncertainZeroParametrization} the main result of Lemma \ref{lem:UncertaintyParameterization} follows, which completes the proof.
\end{proof}

\begin{proof}[\textbf{Proof of Lemma \ref{lem:EstimationErrorBounds}}]
	Since $\Linfnormt{x}{\tau} \leq \rho_x$ and $\Linfnormt{u}{\tau} \leq \rho_u$ by assumption, Lemma \ref{lem:UncertaintyParameterization} holds. Combining \eqref{eq:VirtualSystem} and \eqref{eq:NonlinearParam} yields
	\begin{equation}\label{eq:NonlinearEq1}
	\begin{split}
	\dot{x}_v (t) = & A_m x_v (t) + \bar{B} ( \omega u_z(t)  +  \theta(t) \norm{x_g(t)} + \sigma(t) ), \\
	y(t) = & C_m x_v(t), \quad x_v(0)=x_0  ,
	\end{split}
	\end{equation}
	where $x_g(t) = [x_v^\top(t), x_u^\top(t)]$, and $x_u(t)$, $u_z(t)$ are given in \eqref{eq:ZsU}. Notice that from Corollary \ref{cor:InteactorAugmentedSystem} and Equations \eqref{eq:ZsU}~-~\eqref{eq:barUncertainFunc} it follows that $x(t) = x_v(t) + T_z ( \omega x_u(t) + x_f(t) )$ holds, where $x_f(t)$ is defined in \eqref{eq:ZsF}, and $T_z \in \IR^{n \times n_f}$ satisfies \eqref{eq:InteratorEquiv}.
	Let $v(t) = (\II_n - H C_m) x_v(t)$. By pre-multiplying both sides of \eqref{eq:NonlinearEq1} by $(\II_n - H C_m)$ and taking the derivative of $y(t)$, it follows that
	\begin{align}\label{eq:AugmentedNonlinear}
	\dot{v}(t) = & A_H  v(t) + A_H H y(t) , \quad v(0) = v_0, \nonumber \\
	\dot{y}(t)   = & C_m A_m v(t) + C_m A_m H y(t) \\
	& + C_m \bar{B} ( \omega u_z(t)  + \theta(t) \norm{x_g(t)} + \sigma(t) ), \quad y(0)=y_0 ,\nonumber
	\end{align}
	where $v_0 = (\II_n - H C_m)x_0 $, and $\{A_H, H\}$ is given in \eqref{eq:DefAH}. Let
	\begin{equation}\label{eq:DefEtat}
	\tilde{\eta}_t(t) = \tilde{\omega}(t) u_z(t) + \tilde{\theta}(t) \norm{ \hat{x}_g(t) } + \tilde{\sigma}(t),
	\end{equation}
	and
	\begin{equation}\label{eq:DefPhi}
	\phi(t) = \theta(t) ( \norm{\hat{x}_g(t)} - \norm{{x_g}(t)} ),
	\end{equation}
	where $\tilde{\omega}(t) = \hat{\omega}(t) - \omega$, $\tilde{\theta} (t)= \hat{\theta}(t) - \theta(t)$, and $\tilde{\sigma}(t) = \hat{\sigma}(t) - \sigma(t)$.
	Let
	\begin{equation}\label{eq:DefTildeYV}
	\tilde{v}(t) = \hat{v}(t) - v(t), \quad \tilde{y}(t) = \hat{y}(t) - y(t) .
	\end{equation}
	Then, subtracting \eqref{eq:AugmentedNonlinear} from \eqref{eq:StatePredictor} yields
	\begin{equation}\label{eq:ErrorDynamics}
	\begin{split}
	\dot{\tilde{v}}(t)   = & A_v \tilde{v}(t) -P_v^{-1} A_m^\top C_m^\top P_y \tilde{y}(t) , \\
	\dot{\tilde{y}}(t)   = & -\alpha \tilde{y}(t) + C_m A_m \tilde{v}(t)  + C_m \bar{B} ( \tilde{\eta}_t(t) + \phi(t) )  , \\
	\tilde{v}(0) = & -v_0,    \quad  \tilde{y}(0) = 0,
	\end{split}
	\end{equation}
	where $A_v$ is Hurwitz (see \eqref{eq:DefAv}), and $\tilde{\eta}_t(t)$, $\phi(t)$ are given in \eqref{eq:DefEtat}, and \eqref{eq:DefPhi}, respectively.
	Consider the Lyapunov function:
	\begin{equation}\label{eq:LypunovFunction}
	\begin{split}
	V(t) = & \tilde{v}^\top(t) P_v \tilde{v}(t) + \tilde{y}^\top(t) P_y \tilde{y}(t) \\
	& + \frac{ \tilde{\omega}^2(t)  }{ \Gamma_\omega  } + \frac{ \tilde{\theta}^\top(t) \tilde{\theta}(t)}{ \Gamma_\theta  }
	+  \frac{ \tilde{\sigma}^\top(t) \tilde{\sigma}(t)}{ \Gamma_\sigma  }.
	\end{split}
	\end{equation}
	Taking the derivative of \eqref{eq:LypunovFunction}, and substituting \eqref{eq:AdaptiveLaws} and \eqref{eq:ErrorDynamics}, one has
	\begin{equation}\label{eq:DerivLypunov2}
	\begin{split}
	\dot{V}(t) \leq & - \tilde{v}^\top(t) Q \tilde{v}(t) - 2 \alpha \tilde{y}^\top(t) P_y \tilde{y}(t) - \frac{ 2 \tilde{\theta}^\top(t) \dot{{\theta}} (t)}{ \Gamma_\theta  } \\
	&   - \frac{ 2 \tilde{\sigma}^\top(t) \dot{{\sigma}} (t)}{ \Gamma_\sigma  } + 2\tilde{y}^\top(t) P_y C_m \bar{B} \phi(t),
	\end{split}
	\end{equation}
	where $Q \succ \epsilon_q \II_n $   is positive definite matrix satisfying \eqref{eq:LyapunovEquation}.
	Notice that $\norm{\norm{\hat{x}_g(t)} -  \norm{{x}_g(t)}} \leq \norm{\tilde{v}(t)}$ holds.
	Then, from \eqref{eq:TheSigBound} and \eqref{eq:DefPhi} it follows that
	\begin{equation}\label{eq:PhiBounds}
	\begin{split}
	2\tilde{y}^\top(t) P_y C_m \bar{B} \phi(t) \leq & \alpha_\phi \tilde{y}^\top(t) P_y \tilde{y}(t) \\
	& + \frac{ m \bar{d}_{\rho_x}^2 \twonorm{ \sqrt{P_y} C_m \bar{B}}^2  }{\alpha_\phi} \twonorm{\tilde{v}(t)}^2,
	\end{split}
	\end{equation}
	where $\bar{d}_{\rho_x}$, $\alpha_\phi$ are given in \eqref{eq:DBBound} and \eqref{eq:AlphaY}, respectively.
	Further, from \eqref{eq:DerivLypunov2} and \eqref{eq:PhiBounds} one has
	\begin{equation}\label{eq:DerivLypunov3}
	\begin{split}
	\dot{V}(t) \leq & - \tilde{v}^\top(t) Q_v \tilde{v}(t) - \alpha_y \tilde{y}^\top(t) P_y \tilde{y}(t) \\
	& - \frac{ 2 \tilde{\theta}^\top(t) \dot{{\theta}} (t)}{ \Gamma_\theta  }
	- \frac{ 2 \tilde{\sigma}^\top(t) \dot{{\sigma}} (t)}{ \Gamma_\sigma  } ,
	\end{split}
	\end{equation}
	where $Q_v = Q - \epsilon_q \II_n \succ 0$ and $\alpha_y > 0$ (see \eqref{eq:AlphaY}).
	Notice that from Lemma \ref{lem:UncertaintyParameterization} it follows that for $0 \leq t \leq \tau$
	\begin{equation*}\label{eq:BoundsOnDerivLyapunov}
	\frac{ 2 \tilde{\theta}^\top(t) \dot{{\theta}} (t)}{ \Gamma_\theta  } + \frac{ 2 \tilde{\sigma}^\top(t) \dot{{\sigma}} (t)}{ \Gamma_\sigma  } \leq
	\frac{\theta_1 - \theta_0}{\Gamma} \lambda_1,
	\end{equation*}
	and the projection operator in \eqref{eq:AdaptiveLaws} ensures
	\begin{equation}\label{eq:BoundsOnLyapunov}
	\frac{ \tilde{\omega}^2(t)  }{ \Gamma_\omega  } + \frac{ \tilde{\theta}^\top(t) \tilde{\theta}(t)}{ \Gamma_\theta  }
	+  \frac{ \tilde{\sigma}^\top(t) \tilde{\sigma}(t)}{ \Gamma_\sigma  } \leq \frac{ \theta_0 }{\Gamma} ,
	\end{equation}
	where $\Gamma$, $\theta_0$, $\theta_1$, $\lambda_1$ are given in \eqref{eq:EstimationVariables}.
	Since
	\begin{equation*}\label{eq:InequilityQuadartic}
	\begin{split}
	- \tilde{v}^\top(t) Q_v & \tilde{v} (t) - \alpha_y  \tilde{y}^\top(t) P_y \tilde{y} (t) \\
	& \leq  - \lambda_1 ( \tilde{v}^\top(t) P_v \tilde{v} (t) + \tilde{y}^\top(t) P_y \tilde{y} (t) ),
	\end{split}
	\end{equation*}
	combining \eqref{eq:DerivLypunov3}~-~\eqref{eq:BoundsOnLyapunov}, along with \eqref{eq:LypunovFunction}, leads to
	\begin{equation*}
	\dot{V}(t) \leq -\lambda_1 \left(V(t) - \frac{\theta_1}{\Gamma} \right).
	\end{equation*}
	Choose $t_0 \in \IR$ to be $0 \leq t_0 \leq t \leq \tau$. Then, Gronwell-Bellman inequality yields
	\begin{equation}\label{eq:IneVtandVarepsilon}
	\sqrt{V(t)} \leq \upsilon_v(t, t_0), \quad 0 \leq t_0 \leq t \leq \tau ,
	\end{equation}
	which gives
	\begin{equation}\label{eq:IneLapunove}
	\norm{\tilde{v}(t)} \leq \frac{\upsilon_v(t, t_0)}{\sqrt{\lambda_{\min}(P_v)}}, \quad \norm{\tilde{y}(t)} \leq \frac{\upsilon_v(t, t_0)}{\sqrt{\lambda_{\min}(P_y)}} ,
	\end{equation}
	where
	\begin{equation}\label{eq:DefUpsionV}
	\upsilon_v(t, t_0) = \sqrt{ \left( V(t_0) - \frac{\theta_1}{\Gamma } \right)e^{-\lambda_1(t-t_0) }+ \frac{\theta_1}{\Gamma} } .
	\end{equation}
	Finally, since $V(0) \leq x_0^\top \bar{P}_v x_0 + \frac{\theta_0}{\Gamma}$ with $\bar{P}_v = (\II_n - HC_m)^\top P_v (\II_n - HC_m)$, from letting $t_0 =0$ it follows that
	\begin{equation}\label{eq:EstimationErrorBoundTotal}
	\begin{split}
	\norm{\tilde{y}(t)} \leq & \kappa_y e^{-\frac{\lambda_1}{2} t }\norm{x_0} + \sqrt{ \frac{ \theta_1  }{ \lambda_{\min}(P_y)}  }  \frac{1}{\sqrt\Gamma} , \\
	\norm{\tilde{v}(t)} \leq & \kappa_v e^{-\frac{\lambda_1}{2} t }\norm{x_0} + \sqrt{ \frac{ \theta_1  }{ \lambda_{\min}(P_v)}  }  \frac{1}{\sqrt\Gamma} ,
	\end{split}
	\end{equation}
	where $\kappa_v$, $\kappa_y$ are given in \eqref{eq:DefKappa}.
	This completes the proof.
\end{proof}

\begin{proof}[\textbf{Proof of Theorem \ref{thm:ClosedLoopStability}}]
	Let $\txref (t)= x_{ref}(t) -x(t) $, $\turef(t) = \uref(t) -u(t)$, $ \tyref (t) = \yref(t) -y(t)$, and $\tnref(t) = f(\xref,t) - f(x,t)$. First, it will be shown that Equation \eqref{eq:TransientBounds1} holds by a contradiction argument. Suppose it is not true. Notice that since $\kappa_m \geq 1$ in \eqref{eq:DefKappa}, it follows that $\gamma_{x_0} > 1$, which leads to $\rho_{dx} > \rho_0$, and $ \norm{\txref(0)} = \rho_0 < \rho_{dx}$, where $\gamma_{x_0}$, $\rho_{dx}$ are given in \eqref{eq:GammaXU}, and \eqref{eq:RhoDXU}, respectively. Moreover, since $\norm{\turef(0)} = 0 < \rho_{du} $ with $\rho_{du}$ being given in \eqref{eq:RhoDXU}, from the continuity of the solutions it follows that there exists $\tau' >0$ such that
	\begin{equation*}\label{eq:tXrefAssumption1}
	\norm{\txref(\tau')} = \rho_{dx}\quad \text{or} \quad \norm{\turef(\tau')} = \rho_{du},
	\end{equation*}
	while $\norm{\txref(t)} < \rho_{dx}$ and $\norm{\turef(t)} < \rho_{du}$ for $0 \leq t < \tau'$.
	This implies that the following must hold:
	\begin{equation}\label{eq:tXrefAssumption3}
	\Linfnorm{\txrefsub{\tau'}} \leq \rho_{dx} , \quad \Linfnorm{\turefsub{\tau'}} \leq \rho_{du}.
	\end{equation}
	Notice that from \eqref{eq:DefRhox}, \eqref{eq:DefRhoXR}, and \eqref{eq:RhoDXU} it follows that
	\begin{equation*}
	\rho_{dx} = \rho_x - \rho_{rx} , \quad \rho_{du} = \rho_{u} - {\rho}_{ru}.
	\end{equation*}
	Then, the triangular inequalities on \eqref{eq:tXrefAssumption3}, together with \eqref{eq:RefBoundedX} and \eqref{eq:RefBoundedU}, yield
	\begin{equation}\label{eq:XUbound}
	\Linfnormt{x}{\tau'}  \leq \rho_x , \quad \Linfnormt{u}{\tau'} \leq \rho_u   ,
	\end{equation}
	which, together with Assumption \ref{as:UncertaintiesOfFunction} and the fact that $d_{\rho_x} \leq L_{\rho_r}$, lead to
	\begin{equation}\label{eq:NrefLipschitz}
	\norm{\tnref (t) } \leq L_{\rho_r} \norm{\txref (t)} , \quad 0 \leq t \leq \tau' .
	\end{equation}
	Since Equation \eqref{eq:XUbound} holds, from Lemma \ref{lem:UncertaintyParameterization}, Equation \eqref{eq:DefHatEtaT} can be rewritten as
	\begin{equation}\label{eq:DefEtaHat1}
	{\hat \eta _t}(t) = \omega u_z(t) + \bar{f}(X, t)  + \tilde{\eta}_t(t) + \phi(t) ,
	\end{equation}
	where $u_z(t)$, $\bar{f}(X, t)$, $\tilde{\eta}_t(t)$, $\phi(t)$ are given in \eqref{eq:ZsU}, \eqref{eq:barUncertainFunc}, \eqref{eq:DefEtat}, and \eqref{eq:DefPhi}, respectively.
	Notice that $x(t) = T_gx_g(t) + T_z x_f(t)$, and therefore from \eqref{eq:ZsU} and \eqref{eq:ZsF} it follows that
	\begin{equation}\label{eq:DefEtaHat2}
	\omega u_z(t) + \bar{\eta}(s)= Z(s) ( \omega u(s) + \eta(s) ),
	\end{equation}
	where $\bar{\eta}(s)$, $\eta(s)$ are the Laplace transforms of $\bar{f}(X, t)$ and $f(x,t)$, respectively. Now, substituting \eqref{eq:DefEtaHat1} and \eqref{eq:DefEtaHat2} into \eqref{eq:ControlLaws} leads to
	\begin{equation}\label{eq:L1ControlLawWithEtat}
	\begin{split}
	u(s) = & C_0(s) ( K_g r(s) - \eta (s))  - \phi_c(s), \\
	\phi_c(s) = & C_0(s) Z^{-1}(s) ( {\tilde \eta _t}(s) + \phi(s) ) ,
	\end{split}
	\end{equation}
	where $C_0(s)$ is given in \eqref{eq:DefC0}; $C_0(s) Z^{-1}(s)$ is a stable and strictly proper transfer matrix. Combining the Laplace transform of \eqref{eq:nonlinearmodel1} with \eqref{eq:L1ControlLawWithEtat}
	yields
	\begin{align}\label{eq:DefClosedLoop}
	x(s)  = &  H_r(s)r(s) + G(s)\eta (s) - \omega H_0(s) \phi_c(s) + x_{in}(s), \nonumber \\
	y(s) = & C_m x(s),
	\end{align}	
	where $H_r(s)$, $H_0(s)$, $G(s)$ are given in \eqref{eq:DefHsfun}, and $x_{in}(s) = (s\II_n-A_m)^{-1}x_0 $.
	By subtracting \eqref{eq:L1ControlLawWithEtat} and \eqref{eq:DefClosedLoop} from \eqref{eq:TfReferenceSystem}, it follows that
	\begin{equation}\label{eq:DifSysTf}
	\begin{split}
	\txref(s) = & G(s) \tnref(s) + \omega H_0(s) \phi_c(s) - x_{in}(s), \\
	\tyref(s) = & C_m \txref(s) ,
	\end{split}
	\end{equation}
	and
	\begin{equation}\label{eq:DifInputTf}
	\turef(s) =  - {C_0}(s)\tilde{\eta}_{ref}(s) + \phi_c(s)  .
	\end{equation}
	Since $(C_mB_m)^\dag (C_mB_m) = \II_m$, from \eqref{eq:ErrorDynamics} one has
	\begin{equation}\label{eq:tfCsEquation}
	\begin{split}
	\phi_c(s) =  C_1(s) \ytilde(s) - C_2(s) \vtilde(s),
	\end{split}
	\end{equation}
	where $\{C_1(s), C_2(s)\}$, and $\set{\tilde{y}(t) ,\tilde{v}(t) }$ are defined in  \eqref{eq:DefCsfun1}, and \eqref{eq:DefTildeYV}, respectively; $C_1(s)$, $C_2(s)$ are all stable and proper transfer function matrices.
	From \eqref{eq:L1Condition} it can be shown that $\Lonenorm{G(s)}L_{\rho_r} < 1$. Therefore, combining \eqref{eq:NrefLipschitz}, and \eqref{eq:DifSysTf}-\eqref{eq:tfCsEquation} yields
	\begin{align*}\label{eq:DifSysIneq1}
	\Linfnorm{\txrefsub{\tau'}} \leq & \frac{ \Lonenorm{H_1(s) }\Linfnormt{\tilde{y}}{\tau'}  + \Lonenorm{H_2(s) }\Linfnormt{\tilde{v}}{\tau'}  }{1 - \Lonenorm{G(s)} L_{\rho_r} }  \nonumber \\
	& + \frac{ \kappa_m \rho_0 }{1 - \Lonenorm{G(s)} L_{\rho_r}} , \\
	\Linfnorm{\turefsub{\tau'}} \leq &   \Lonenorm{C_1(s)  }\Linfnormt{\tilde{y}}{\tau'} + \Lonenorm{C_2 (s) }\Linfnormt{\tilde{v}}{\tau'} \nonumber \\
	&  + \Lonenorm{C_0(s) } L_{\rho_r} \Linfnorm{\txrefsub{\tau'}}  ,\nonumber
	\end{align*}
	where $\kappa_m$, $\set{H_1(s), H_2(s)}$, $\set{C_1(s), C_2(s)}$ are given in \eqref{eq:DefKappa}, \eqref{eq:DefHsfun}, and \eqref{eq:DefCsfun1}, respectively. Since Equation \eqref{eq:EstimationErrorBoundTotal} holds for $0 \leq t \leq  \tau'$, one has
	\begin{equation}\label{eq:DifInputIneq1}
	\begin{split}
	\Linfnorm{ \turefsub{\tau'} } \leq & \gamma_{u_0} \rho_0 + \gamma_u \sqrt{ \frac{\theta_1}{\Gamma}}, \\
	\Linfnorm{ \txrefsub{\tau'} } \leq & \gamma_{x_0} \rho_0 + \gamma_x \sqrt{ \frac{\theta_1}{\Gamma}} ,
	\end{split}
	\end{equation}
	where $\gamma_{u_0}$, $\gamma_{x_0}$, $\gamma_{u}$, $\gamma_{x}$ are given in \eqref{eq:GammaXU}, and $\theta_1$ is defined in \eqref{eq:EstimationVariables}.
	Since $\Gamma > 0$ is chosen so that $  \gamma_x \sqrt{ \frac{\theta_1}{\Gamma}} < \bar{\gamma}$ and $  \gamma_u \sqrt{ \frac{\theta_1}{\Gamma}} < \bar{\gamma}$, from \eqref{eq:DifInputIneq1} it follows that
	\begin{equation*}\label{eq:DifInputIneq2}
	\Linfnorm{ \turefsub{\tau'} } <  \rho_{dx} , \quad \Linfnorm{ \txrefsub{\tau'} } < \rho_{du},
	\end{equation*}
	which contradict  \eqref{eq:tXrefAssumption3}, thus proving \eqref{eq:TransientBounds1}.
	Moreover, Equation \eqref{eq:TransientBounds2} is obtained from applying the triangular inequality on $\norm{\txref}$ and $\norm{\turef}$.
	
	Nest we prove Equation \eqref{eq:TransientBounds3}.
	Let $A_b \in \IR^{n_b \times n_b}$, $B_b \in \IR^{ n_b \times m}$, and $C_b \in \IR^{ m \times m}$ be a minimal realization of $C(s)$ with the appropriate dimension $n_b$. Then, the system given in \eqref{eq:DifSysTf} and \eqref{eq:DifInputTf} can be represented as
	\begin{equation} \label{eq:AcBcCc} 
	\begin{split}
	\dot{\tilde{x}}_c(t) & = A_c\tilde{x}_c(t) + B_c  \tnref (t)  + \bar{B}_c \omega \phi_c(t) , \\
	\txref(t) & = C_c \tilde{x}_c(t) ,  \quad \tilde{x}_c(0) =[-x_0^\top, 0]^\top,
	\end{split}
	\end{equation}
	with
	\begin{equation*}\label{eq:DefAcBcCc}
	\begin{split}
	{A_c}  = & \left[ {\begin{array}{*{20}{c}}
		{{A_m}}&{{B_m}{C_b}}\\
		0&{{A_b}}
		\end{array}} \right], \quad
	{B_c}  = \left[ {\begin{array}{*{20}{c}}	{{B_m}}\\{{-B_b}}\end{array}} \right], \\
	{\bar{B}_c} = & \left[ {\begin{array}{*{20}{c}}	{{B_m }}\\{{0}}\end{array}} \right], \quad
	{C_c}  =  \left[ {\begin{array}{*{20}{c}}\II_n &0\end{array}} \right] ,
	\end{split}	
	\end{equation*}
	where $\tilde{x}_c(t) = [\txref^\top(t), \tilde{x}_b^\top(t)]^\top \in \IR^{n_c \times n_c}$ is the state vector with $n_c = n+n_b$. Let $t_m \geq 0$. Then, from \eqref{eq:AcBcCc} it follows that for $t \geq t_m$
	\begin{equation}\label{eq:XcSoln}
	\begin{split}
	\tilde{x}_c(t) = & e^{A_c (t-t_m)} \tilde{x}_c(t_m) + \int_{t_m}^{t} {e^{A_c (t-\tau) } \bar{B}_c   \omega \phi_c(\tau)  d\tau } \\
	& + \int_{t_m}^{t} {e^{A_c (t-\tau) } B_c  \tnref(\tau)   d\tau } 	.
	\end{split}	
	\end{equation}
	Notice that it can be shown that $G(s) = C_c ( s\II_{n_c} - A_c)^{-1} B_c = H_0(s)(\II_m - C(s))$. Since $\Lonenorm{G(s)} L_{\rho_r} < 1$ holds from \eqref{eq:L1Condition}, from the continuity of the $\Lone$-norm, one may take a sufficiently small $\lambda_0 >0$ such that $\beta_1 = \Lonenorm{ G(s-\lambda_0) } < 1/L_{\rho_r}$. Let $A_{\lambda_0} = A_c + \lambda_0 \II_{n_c}$, and define $\bar{x}_{c}(t) = e^{\lambda_0 (t-t_m)} \tilde{x}_c(t)$, $\bar{\phi}_c(t) = \omega e^{\lambda_0 (t-t_m)} \phi_c(t)$, $\bar{x}_{ref}(t) = e^{\lambda_0 (t-t_m)} \txref(t)$, and $\bar{\eta}_{ref}(t) = e^{\lambda_0 (t-t_m)} \tnref(t)$.
	Since Assumption \ref{as:UncertaintiesOfFunction} implies that
	\begin{equation}\label{eq:BarEtaRefIneq}
	\norm{ \bar{\eta}_{ref}(t)   } \leq L_{\rho_r} \norm{ \bar{x}_{ref} (t) } ,
	\end{equation}
	multiplying both sides of \eqref{eq:XcSoln} by $e^{\lambda_0 (t-t_m)}C_c$ leads to
	\begin{equation}\label{eq:XrefBar}
	\begin{split}
	\Linfnormtb{\bar{x}_{ref}}{t_m}{t} \leq & \frac{\beta_0 }{1 - \beta_1 L_{\rho_r} }\norm{\tilde{x}_c(t_m)} \\
	& +  \frac{\beta_2}{1 - \beta_1  L_{\rho_r} } \Linfnormtb{ \bar{\phi}_c}{t_m}{t} ,
	\end{split}
	\end{equation}
	where $\beta_0 = \sup\limits_{0 \leq \tau } { \norm{e^{A_{\lambda_0} \tau } } }$, $\beta_1 = \Lonenorm{G(s-\lambda_0)}$, and $\beta_2 = \Lonenorm{  (s\II_{n_c} - A_{\lambda_0})^{-1} \bar{B}_c}$. By combining \eqref{eq:XcSoln}~-~\eqref{eq:XrefBar}, it can be shown that	
	$\norm{\bar{x}_c(t)} \leq  \kappa_0 \norm{\tilde{x}_c(t_m)}
	+ \kappa_1 \Linfnormtb{\bar{\phi}_{c}}{t_m}{t}  $,
	which further gives
	\begin{equation}\label{eq:XcBdd2}
	\begin{split}
	\norm{\tilde{x}_c(t)} \leq & \kappa_0 e^{-\lambda_0 (t-t_m)} \norm{\tilde{x}_c(t_m)} +  \omega_u \kappa_1 \Linfnormtb{{\phi}_{c}}{t_m}{t}  ,
	\end{split}
	\end{equation}		
	where $\omega_u >0$ is the upper bound of $\omega$, and
	\begin{equation}\label{eq:Kappas}
	\begin{split}
	\kappa_0 = & \beta_0 ( 1 + \frac{ L_{\rho_r} \beta_3 }{1 - \beta_1 L_{\rho_r} }) , \quad
	\kappa_1 =   \beta_2 ( 1 + \frac{ L_{\rho_r}  \beta_3 }{1 - \beta_1 L_{\rho_r} }  ),
	\end{split}
	\end{equation}
	with $\beta_3 = \Lonenorm{ (s\II_{n_c} - A_{\lambda_0})^{-1} B_c }$. Substituting \eqref{eq:tfCsEquation}, together with \eqref{eq:IneVtandVarepsilon}~-~\eqref{eq:DefUpsionV}, into \eqref{eq:XcBdd2} leads to
	\begin{equation}\label{eq:XcBdd3}
	\begin{split}
	\norm{\tilde{x}_c(t)} \leq & \kappa_0 e^{-\lambda_0 (t-t_m)} \norm{\tilde{x}_c(t_m)} \\
	&  +  \gamma_1  \left( \upsilon_v(t_m,0) + \sqrt{\frac{\theta_1}{\Gamma}}\right) ,
	\end{split}
	\end{equation}
	where $\upsilon_v(\cdot,\cdot)$, $\theta_1$ are defined in \eqref{eq:DefUpsionV}, and \eqref{eq:EstimationVariables}, respectively,  and
	\begin{equation*}
	\gamma_1 = \omega_u \kappa_1 \left(  \frac{\Lonenorm{C_1(s)}}{ \sqrt{ \lambda_{\min}(P_y) }  } + \frac{\Lonenorm{C_2(s)}}{ \sqrt{ \lambda_{\min}(P_v) }  } \right).
	\end{equation*}
	Notice that from \eqref{eq:EstimationVariables} and \eqref{eq:DefUpsionV} it follows that
	\begin{equation*}
	\upsilon(t_m, 0) \leq \sqrt{ n \lambda_{\max}(\bar{P}_v)} \norm{x_0} + \sqrt{\frac{\theta_1}{\Gamma}},
	\end{equation*}
	which, together with \eqref{eq:XcBdd3}, results in
	\begin{align}\label{eq:XcBdd4}
	\norm{\tilde{x}_c(t_m)} \leq & \kappa_0 e^{-\lambda_0 t_m } \norm{x_0} + \gamma_1 \sqrt{ n \lambda_{\max}(\bar{P}_v) }  \norm{x_0} \nonumber \\ & + \gamma_1 \sqrt{ \frac{\theta_1}{\Gamma}} ,
	\end{align}
	where $\bar{P}_v = (\II_n - H C_m)^\top P_v (\II_n - H C_m)$.		
	Let $t_m=t/2$. Then, substituting \eqref{eq:XcBdd4} into \eqref{eq:XcBdd3}, and using \eqref{eq:DefUpsionV}, one has
	\begin{equation*}
	\norm{\tilde{x}_{ref}(t)} \leq  \norm{\tilde{x}_c(t)} \leq \upsilon_{dx}(t) \norm{x_0} + \frac{\gamma_{dx}  }{\sqrt{\Gamma}} , \quad t \geq 0,
	\end{equation*}	
	where
	\begin{equation*}
	\begin{split}
	\upsilon_{dx} (t) = & \kappa_0^2  e^{ - \lambda_0 t } + \kappa_0 \gamma_1 \sqrt{ n \lambda_{\max}(\bar{P}_v) } e^{ -\frac{\lambda_0 }{2} t} \\
	& 	+ \gamma_1 \sqrt{ n \lambda_{\max}(\bar{P}_v) } e^{ - \frac{ \lambda_1 }{4} t }, \\
	\gamma_{dx}    = & ( \kappa_0 + 2 ) \gamma_1 \sqrt{  \theta_1  },
	\end{split}
	\end{equation*}
	with $\lambda_1$ being given in \eqref{eq:EstimationVariables}.
	Finally, letting $\gamma_{dy} = \norm{C_m} \gamma_{dx} $,  $\upsilon_{dy}(t) = \norm{C_m} \upsilon_{dx}(t) $ reduces to \eqref{eq:TransientBounds3}. This completes the proof.	
\end{proof}

\vfill

\end{document}